\documentclass[preprint,12pt]{elsarticle}

\newtheorem{proposicao}{{\bf Proposition}}

\usepackage{amssymb}

\journal{Some Journal in the World}

\begin{document}

\begin{frontmatter}

\title{Mathematical Properties of  Strategies to Control Epidemic Outbreaks in the Context of SEIR Models with Multiple Infectious Stages}

\author[rvt,focal]{Annibal Figueiredo}
\author[rvt,focal]{Tarcísio Marciano da Rocha Filho}

\address[rvt]{Institute of Physics, University of Brasilia, 70919-970 - Brasilia, Brazil.}
\address[focal]{International Center of Condensed Mater, 70919-970 - Brasilia, Brazil}

\begin{abstract}
In this work we analyze mathematically the consequences and effectiveness of strategies to control an epidemic in the framework of classical SEIR models with multiple parallel infectious stages. We define the mathematical concept of a control strategy, showing that it implies turning classic epidemiological models into systems of non-autonomous differential equations.
The analysis of these non-autonomous systems is based on the two main results obtained  in this work: the first establishes a condition that implies a dynamic without epidemic outbreaks; the second establishes a maximum value for the susceptible population associated to the fixed points that are attractors, moreover, we proof that any trajectory converges to some of these attractors.
An important consequence of this last result  is the existence of an insurmountable limit on the number of infected individuals after the end of a given control strategy.
This restriction can only be mitigated by changing the maximum value of susceptible population associated to the system attractors, which could only be done with  permanent control action, that is, without returning to normality.
Another interesting result of our work is to show how the moment to start and the way how the control strategy ends  strongly impacts the asymptotic value for the total number of infected individuals. We illustrate our analysis and results in a SEIR model (with two or three parallel stages) applied to describe the COVID-19 epidemic.
\end{abstract}


\end{frontmatter}

\newpage

\section{Introduction}

The emergence of COVID-19 epidemic in Wuhan (Hubei - China) \cite{ref1,ref2} and its rapid spread to all regions of the world \cite{ref3,ref4}, reaching pandemic status, have put enormous pressure on all healthcare systems in the most diverse countries, including Brazil \cite{ref5}, mainly due to the disease severity and absence of effective pharmaceutical treatment. Therefore, only non pharmaceutical interventions, that relive the classic and secular practices of social isolation and quarantine, remained as strategies to control the epidemic \cite{ref6,ref7}. Many works, produced at a frantic pace, have shown the effectiveness and quality of classic epidemiological models \cite{ref8,ref9,ref10,ref11,ref12} for the assessing and monitoring of epidemic evolution in different regions of the world. This monitoring and anticipation work on how the epidemic evolve is fundamental to support governmental decisions to mitigate the spread and effects of COVID-19 epidemic. The existence of  abundant data and efficient statistical techniques allows the estimation of the fundamental  parameters used in these models with high precision and confidence \cite{ref13,ref14,ref15,ref16,ref17,ref18,ref19}. 

The vast majority of  models used to study and analyze strategies of epidemic control, whether in the deterministic \cite{ref18,ref19,ref20,ref21,ref22,ref23,ref24,ref25,ref26,ref27,ref28,ref29} or stochastic \cite{ref30,ref31,ref32} version, are modifications of the classical SIR or SEIR models with multiple and  infectious stages \cite{ref33,ref33a}: in these models, the infected (SIR model) or exposed (SEIR model) individuals evolve in to different infectious stages. These stages may be related to different characteristics associated to infected individuals such as symptoms, the ability to transmit the disease or, concomitantly and non-exclusively, the greater or lesser capacity for detecting them. In the case of COVID-19 these characteristics were evidenced due the importance of taking into account asymptomatic transmitters or, more precisely, individuals with symptoms that do not allow to distinguish COVID-19 from other typical respiratory infections except through the performance of specific laboratory tests \cite{ref15,ref28,ref34,ref35}.

A deterministic mathematical model of a control strategy of an epidemic is set up through the following steps.
\begin{enumerate}
\item The choice of an autonomous  system of differential equations to model the natural evolution of an epidemic. By the expression ``natural evolution'' we understand an epidemic evolution where no intervention or control measures are made in order to mitigate its propagation. We call this mathematical model as an Epidemic Natural Model (ENM) 
\item The decision in how the chosen control strategy  must be modeled:
\begin{enumerate} 
\item trough possible value modification for some of the parameters appearing in the ENM and, in this case, these parameters are modeled as time dependent functions and the respective ENM becomes a non-autonomous system of differential equations (see references \cite{ref18,ref19,ref20,ref21,ref22,ref23,ref24,ref25}),
\item or trough reformulation of the ENM with incorporation of new variables (and the respective differential equations)  that are associated to the control strategy, leading to a new system of differential equations, where one may also consider that some parameters are time dependent functions (see references \cite{ref18,ref26,ref27,ref28,ref29}).       
\end{enumerate}     
\end{enumerate}

Control models in (a) are generally related to measures and actions that impact the basic reproduction rates of the different infectious stages. When this impact is the same for all infectious stages, the strategy is said to be uniform. As examples of this kind of strategies, we can consider the use of masks by all the population and strategies of social distancing: from the mildest, as to avoid agglomerations or maintain a minimum distance between individuals, even the most rigorous ones, such as the lockdown.
On the other hand, as example of a control model in (b) we have the idea of quarantine, which involves identifying and isolating infected individuals. In fact, this kind of control strategy is the standard and traditional way of dealing with infectious diseases and it is linked to the greater or lesser capacity of epidemiological surveillance of public healthcare systems. For  COVID-19, this capacity of epidemiological surveillance has been under high pressure due the difficulty to detect and monitor infected individuals with mild or no symptoms of the disease.
However, it is worth to emphasize that the same kind of control strategy has not necessarily just one way to be mathematically modeled.

The main objective of this work is to develop some critical reflections about  epidemic control strategies that are based on deterministic mathematical models. To conduct such study, we develop some general ideas about control strategies in the context of a SEIR model with multiple infectious stages, as originally defined in \cite{ref33}, and we establish some general properties. Amid these properties, we emphasize the two main  results obtained and demonstrated in this work: the first establishes a condition that implies a dynamic without epidemic outbreaks; the second establishes a maximum value for the susceptible population associated to the fixed points that are attractors, moreover, we proof that any trajectory, associated to a given initial condition, converges to some of these attractors. This last result allows to formulate what we call  ''Fundamental Principle of  Finite Control Strategies'',  which affirms the existence of a maximum value (strictly less than one) for the asymptotic limit of susceptible individuals.

An important consequence of this Fundamental Principle  is to show the existence of a minimum limit for the total final number of infected individuals after the return to normality, that is,
the asymptotic value of the total number of infected individuals after the end of the control. Indeed, this demonstrates an insurmountable limit about the minimum number of infected individuals after the end of control strategies.
The restriction imposed by this result can only be mitigated by changing the maximum limit of  susceptible associated to the system attractors, which could only be done with  permanent control action, that is, without returning to normality. The only viable way to do this is through the permanent maintenance of a good epidemiological surveillance system, otherwise we would have to permanently reduce the basic reproduction numbers of the different infectious stages with measures that impact the entire population.

The study of control strategies seek to focus on three factors are fundamental to analyze the effectiveness of a control strategy: trying to minimize the total number of infected individuals after the end of the epidemic cycle, avoiding peak and very high outbreaks that compromise the treatment capacity of all ill individuals and finally the duration of  control strategy, which always implies high social and economic costs. The extensions of the analysis to a SIR model with multiple infectious stage is straightforward and are not present here in this article.
We illustrate these ideas in a model for COVID-19 with two infectious stages (SEIAR model without isolation of infected individuals) and a model with three infectious stages (SEIAQR model with isolation of some of the infected individuals- quarantine).

Another interesting result of our work is to show how the moment to start and the way how the control strategy ends  strongly impacts the asymptotic value for the total number of infected individuals. We show this result through the analysis of finite strategies that diminishes the basic reproduction numbers over a period of time.
We also show how finite strategies, that reduce the basic reproduction numbers to a certain fixed values, imply the existence of a new peak after the end of the control.
To avoid the existence of a second peak, we defined what we call non-outbreak and weak outbreak strategies that prevent the existence of strong outbreak after the controlling finishes.

Finally, it follows from the results obtained in our work that it is not unequivocally defined the value of the so called herd immunity \cite{ref36,ref37} for SEIR models submitted to a finite control strategy. In other words, this means that the total asymptotic number of susceptible individuals, after the end of the control, depends on the control strategy applied. The only thing that can be said is that this value is limited by the fundamental principle of finite strategies.

This article is organized into the following sections: in section 2 we define the SEIR  model with parallel stages and establish some basic properties; in section 3 we define the concept of control strategy and formulate the fundamental principle associated with finite strategies; in section 4 we define a SEIAR model for COVID-19 and analyze the effects of uniform control strategies; in section 5 we define the SEIAQR model for strategies that imply the isolation (quarantine) of infected individuals  and in section 6 we close the article with the conclusions and final considerations. One appendix is devoted to proof the proposition 5 presented in the article.

\section{The SEIR Model with $N$ Multiple and Parallel Infectious Stages} 

One of the most interesting aspects in mathematical modeling of  COVID-19 epidemic  using a typical SEIR system is the need to separate the infected population into two different stages: the symptomatic infected and the asymptomatic infected. This fact leads to the formulation of a SEIAR model that divides the variable $ I $, associated with the infected population, in two variables: one continues to be denoted by $ I $, but now it represents the symptomatic infected, and the other variable, denoted by $ A $, represents the asymptomatic infected \cite{ref42}.

The SEIAR model is just a very particular case of a more general model in which the infected population is divided into $ N $ different stages of infected individuals, which differ from each other only because they have different basic reproduction numbers. We define this generalized model through the following system of differential equations \cite{ref43}:
\begin{equation}
\label{seir}
\begin{array}{c}
 {\dot S}=-\left(\sum_{i=1}^{N}\beta_iI_i\right)S,\;\;{\dot E}=\left(\sum_{i=1}^{N}\beta_iI_i\right)S-\sigma E,\\ \\
 {\dot I}_i=x_i\sigma E-\gamma_iI_i,\;\;{\dot R}=\sum_{i=1}^N\gamma_iI_i,
\end{array}\;(i=1,\ldots,N),
\end{equation}
where the meaning of variables and parameters appearing in the system constituted by $N+3$ ODE's is summarized in table 1.

\vspace*{0mm}
\begin{center}
\begin{tabular}{c}
Table 1 - Variables and Parameters in ODE's system (\ref{seir}) \\
\hline\hline 
Model Variables ($i=1,\ldots,N$)\\
\hline\hline 
\begin{tabular}{cc}
$S$ & Susceptible Population\\ \\
$E$ & Exposed Population \\ \\
$I_i$ & Infected Population in Stage $i$ \\ \\
$R$ & Recovered Population
\end{tabular} \\ 
\hline\hline
Model Parameters - $\beta_i=\gamma_i R_i$ ($i=1,\ldots,N$) \\ \hline\hline
\begin{tabular}{cc}
$R_i$ & Basic Reproduction Number of  Infected in Cohort $i$ \\  \\
$\sigma^{-1}$ & Incubation Time  \\ \\
$\gamma_i^{-1}$ & Infection Time of Infected  Population in Stage $i$ \\ \\
$x_i$ & Likelihood of Exposed becomes an Infected\\
&  in Stage $i$ $\left(\sum_{i=1}^Nx_i=1\right)$ \\
\hline\hline
\end{tabular}
\end{tabular}
\end{center}
The two-dimensional surface $ E = I_1 =\cdots=I_N = 0 $ is an invariant surface of the system (\ref{seir}) and it  is constituted by an infinity number of fixed points, defined by the different possible values of $0\leq S\leq 1$ and $0\leq R\leq 1$. Indeed, this surface of fixed points corresponds to a two-dimensional hyperplane in the $N+3$ dimensional space
$I\hspace*{-1.5mm}R^{N+3}$.
For this kind of  invariant surface, if we define the column vector $ {\vec V} = [E, I_1,\ldots,I_N] ^ {t} $, then, as shown in references  \cite{ref40} and \cite{ref41}, we can write
\begin{equation}
\label{dVdt}
\frac{d{\vec V}}{dt}={\bf L}{\vec V},\;\;{\bf L}(S)=\left(
\begin{array}{ccccc}
-\sigma & \beta_1 S & \beta_2 S & \cdots & \beta_N S\\ & & & &\\
x_1 \sigma & -\gamma_1 & 0  & \cdots & 0 \\ & & & &\\
x_2 \sigma & 0 & -\gamma_2 & \cdots & 0  \\ & & & & \\
\vdots & \vdots  &\vdots &\ddots & \vdots \\ & & & &\\  
x_N\sigma &  0 & 0 &\cdots & -\gamma_N \\
\end{array}\right).
\end{equation}

The main goal of this section is to define some special values for the $S$ variable, which are important to characterize two different dynamical condition in the evolution of the epidemic: the {\it non-outbreak condition} and the {\it weak outbreak condition}. 

\subsection{The non-outbreak condition}

The integral solution of  the linear non-autonomous differential system (\ref{dVdt}), with initial condition ${\vec V}(t_0)$ for some time $t_0>0$, can be formally written as
\begin{equation}
\label{intV}
{\vec V}(t)=\exp\left(\int_{t_0}^t{\bf L}(t')dt'\right){\vec V}(t_0),
\end{equation}
where the elements $ G_ {ij} \neq 0 $ of the matrix $ {\bf G} = \int_ {t_0} ^ t {\bf L} (t ') dt' $ are given by
\begin{equation}
\label{matrizG}
\begin{array}{c}
G_{11}=-\sigma\left(t-t_0\right),\;\;G_{ii}=-\gamma_{i-1}\left(t-t_0\right),\\
\displaystyle G_{1i}=\beta_{i-1}\int_{t_0}^t S(t')dt',\;\;
\displaystyle G_{i1}=x_{i-1}\sigma \left(t-t_0\right),
\end{array}\;\;(i=2,\ldots,N+1).
\end{equation}
It is easy to see that, due to the fact that $ dS/dt '\leq 0 $ in (\ref{seir}), the integrals in the equation (\ref{matrizG}) can be superiorly limited as follows
\begin{equation}
\label{majoracaoS}
 \displaystyle G_{1i}=\beta_{i-1}\int_{t_0}^t S(t')dt'\leq \beta_{i-1}S_0\left(t-t_0\right)\; (i=2,\ldots,N+1),\;S_0=S(t_0).
\end{equation}  

To proceed, let us remember the definition of the norm of the maximum associated to a given matrix:
\begin{equation}
\label{normaM}
|{\bf M}|=\max\left\{|M_{ij}|\right\},\;\;{i=1,\ldots, n,\;j=1\ldots m},
\end{equation} 
where $ M_ {ij} $ are the elements of the $n\times m $ matrix $\bf M $.

From equation (\ref{intV}), the properties of the norm, the relations in (\ref{matrizG}) and the inequality given in (\ref{majoracaoS}), we can show that
\begin{equation}
\label{majoracaoV1}
|{\vec V(t)}|\leq \left|\exp({\bf G})\right| |{\vec V}(t_0)|\leq\left|\exp\left((t-t_0){\bf L}_0\right)\right| |{\vec V}(t_0)|,
\end{equation}
where $ {\bf L} _0 $ is the $ {\bf L} (S) $ matrix defined in (\ref{dVdt}), taking $ S = S_0 $ and $ \gamma_i = \gamma \; (i = 1, \ldots, N) $, where $ \gamma $ is the smallest of the values of $ \gamma_i $. Indeed, we can write explicitly:
\begin{equation}
\label{matrizL0}
{\bf L}_0=\left(
\begin{array}{ccccc}
-\sigma & \beta_1 S_0 & \beta_2 S_0 & \cdots & \beta_N S_0\\ & & & &\\
x_1 \sigma & -\gamma & 0  & \cdots & 0 \\ & & & &\\
x_2 \sigma & 0 & -\gamma & \cdots & 0  \\ & & & & \\
\vdots & \vdots  &\vdots &\ddots & \vdots \\ & & & &\\  
x_N\sigma &  0 & 0 &\cdots & -\gamma \\
\end{array}\right),\;\gamma=\min\{\gamma_i\}_{i=1,\ldots,N}.
\end{equation} 

It is straightforward to show that the matrix $ {\bf L}_0 $ has $ N + 1 $ real eigenvalues written as:
\begin{equation}
\label{eigenL}
\begin{array}{c}
\displaystyle \lambda_1(S_0)= -\frac{\sigma+\gamma}{2}+\frac{1}{2}\sqrt{4\sigma\left(\sum_{i=1}^N\beta_ix_i\right)S_0+(\gamma-\sigma)^2}\\ \\
\displaystyle \lambda_2(S_0)= -\frac{\sigma+\gamma}{2}-\frac{1}{2}\sqrt{4\sigma\left(\sum_{i=1}^N\beta_ix_i\right)S_0+(\gamma-\sigma)^2} \\ \\
\lambda_i=-\gamma\;\;(i=3\cdots N+1).
\end{array}
\end{equation}

If we define $ {\bf D} _0 $ as the diagonal matrix formed by these eigenvalues,
then, from the properties of diagonalizable matrices, it can be shown that there is an invertible matrix $ {\bf T} _0 $ such that
$ {\bf L} _0 = {\bf T} _0 ^ {- 1} {\bf D_0} {\bf T} _0 $,
where the matrix $ {\bf T} _0 $ also depends on $ S_0 $.
The relation between $ {\bf D} _0 $ and $ {\bf L} _0 $ is called similarity relation and implies a similarity relation between the exponentials of $ (t-t_0) {\bf D} _0 $ and $ (t-t_0) {\bf L} _0 $, that is,
\begin{equation}
\label{similarexp}
\exp\left((t-t_0){\bf L}_0\right)={\bf T}_0^{-1}\exp\left((t-t_0){\bf D}_0\right){\bf T}_0.
\end{equation} 
Substituting the relation (\ref{similarexp}) in (\ref{majoracaoV1}) and using the properties of the norm, we can obtain:
\begin{eqnarray}
\label{majoracaoV2}
|{\vec V}(t)|&\leq&\left|{\bf T}_0^{-1}\exp\left((t-t_0){\bf D}_0\right){\bf T}_0\right| |{\vec V}(t_0)|\nonumber \\
&\leq& \left|\exp\left((t-t_0){\bf D}_0\right)\right|
\left|{\bf T}_{0}^{-1}\right|\left|{\bf T}_{0}\right||{\vec V}(t_0)|.
\end{eqnarray}
Noting that the eigenvalue $ \lambda_1 (S_0) $ in (\ref{eigenL}) is the largest eigenvalue of $ {\bf L} _0 $ (respectively of $ {\bf D}_0 $), then the inequality in (\ref{majoracaoV2}) becomes
\begin{equation}
\label{majoracaoV3}
|{\vec V}(t)|\leq \left|{\bf T}_{0}^{-1}\right|\left|{\bf T}_{0}\right||{\vec V}(t_0)| e^{\displaystyle \lambda_1(S_0)(t-t_0)},\; \forall t\geq t_0.
\end{equation}
If we consider the euclidean norm $||\vec V(t)||=\sqrt{E^2(t)+\sum_{i=1}^NI_i^2(t)}$, then the inequality in (\ref{majoracaoV3}) leads to
\begin{equation}
\label{majoracaoV4}
||{\vec V}(t)||\leq 
(N+1)\left|{\bf T}_{0}^{-1}\right|\left|{\bf T}_{0}\right||{\vec V}(t_0)| e^{\displaystyle \lambda_1(S_0)(t-t_0)},\; \forall t\geq t_0.
\end{equation}

It is easy to verify that the eigenvalues of ${\bf L}_0$, with the exception of $ \lambda_1 (S_0) $, are negative for any value of $ S_0 $.
Now, we want to analyze under which conditions $ \lambda_1 (S_0) $ is negative, or, in other words, for which values of $ S_0 $ all eigenvalues $ {\bf L} _0 $ are negative.
It is also easy to verify that the largest eigenvalue $ \lambda_1 (S_0) $ is negative when $ S_0 = 0 $, positive when $ S_0 = 1 $ and changes signs when the determinant of $ {\bf L} _0 $ is null. Hence, we define $ {\bar S} $ as being the value of $ S $ for which  the determinant of $ {\bf L} _0 $ is null, that is,
\begin{equation}
\label{SO}
\det\left({\bf L}_0({\bar S})\right)=0\;\;\iff\;\; {\bar S}=\gamma\left({\displaystyle\sum_{i=1}^N x_i\beta_i}\right)^{-1}.
\end{equation}
Finally, we can conclude that
\begin{equation}
\lambda_1(S_0)<0\;\;\iff\;\;0\leq S_0 < {\bar S}
\end{equation}
and, consequently, from the relation given in (\ref {majoracaoV4}), the value of $ || {\vec V} (t) || $ is limited by an exponential function that is decreasing with time and tends to zero when $ t \rightarrow \infty $. The inequality (\ref{majoracaoV3}) is called the {\it non-outbreak condition} when $ \lambda_1 (S_0) <0 $.

Therefore, since the  $S(t)$ is a non-increasing function, we can establish the following result:

\begin{proposicao}(The Non-Outbreak Condition) The non-outbreak condition is met for every $ t> t_0 $ such that $ S_0 = S (t_0) <{\bar S} $.\end{proposicao}

\subsection{The weak outbreak condition}

Let us consider the function
\begin{equation}
\label{Uliapunov}
U=E+{Z}\sum_{i=1}^N\frac{\beta_i}{\gamma_i}I_i,
\end{equation}
where ${Z}>0$. Tacking the time derivative of $U$ and substituting the expression for $\dot E$ and ${\dot I}_i$ given in (\ref{seir}), we obtain:
\begin{equation}
\label{dUdt}
{\dot U}=(S-{Z})\sum_{i=1}^N\beta_iI_i-\sigma E\left(1-\frac{Z}{{S^*}}\right),
\end{equation}
with ${S^*}$ given by
\begin{equation}
\label{SL}
 {S^*}=\left(\sum_{i=1}^N\frac{x_i\beta_i}{\gamma_i}\right)^{-1}=\left(\sum_{i=1}^N x_i R_i\right)^{-1},
\end{equation} 
where, to obtain the last expression for ${S^*}$, we have used $\beta_i=\gamma_iR_i$ as defined in table 1. Let us observe that the value of ${S^*}$ corresponds to the value of $S$ for which the determinant of the matrix ${\bf L}(S)$ in (\ref{dVdt}) is null. 

If we take ${Z}={S^*}$, then the expression for ${\dot U}$ in (\ref{dUdt}) becomes
\begin{equation}
\label{dUdt1}
{\dot U}=(S-{S^*})\sum_{i=1}^N\beta_iI_i.
\end{equation}
Now, let us suppose that $U(0)>0$ for an initial instant of time $t=0$ such that $S(0)>{S^*}$. As the function $S(t)$ is decreasing and $S(t)\geq 0$ for all $t>0$,  we have that 
$0\leq S_\infty=\lim_{t\rightarrow\infty}S(t)$ exists and $\lim_{t\rightarrow\infty}{\dot S}(t)=0$.  

Also, let us suppose that $S_\infty\geq {S^*}$. Therefore, we can conclude that ${\dot U}(t)\geq 0$ for all time $t>0$. As $U(t)$ is bounded, then it converges to some positive value.  From the equation for ${\dot S}$ in  (\ref{seir})  we  can conclude that $\lim_{t\rightarrow\infty}I_i(t)=0$ for all $i=1,\ldots,N$. Finally, we conclude that $\lim_{t\rightarrow\infty}E(t)=E_\infty>0$.  

Hence, if we suppose that $S_\infty \geq {S^*}$ we should conclude that the considered solution of the system (\ref{seir}) converges to a fixed point where $I_i=0$ and $E=E_\infty>0$. However, this fixed point does not exist  because there are only free disease fixed points associated to the model given by the set of differential equations in (\ref{seir}).     
Thus, we must conclude that $S_\infty <{S^*}$.

The arguments presented in the precedent paragraphs can be synthesized in the following proposition: 

\begin{proposicao}(Maximum Asymptotic Limit for $ S $) For any solution of the system (\ref{seir}), with initial condition in time $ t_0 $ such that $   U (t_0)> 0 $ and $ S (t_0) > {S^*} $, we have
$$S_\infty=\lim_{t\rightarrow\infty}S(t)< {S^*},\;\;\lim_{t\rightarrow\infty}E(t)=0,\;\;\lim_{t\rightarrow\infty}I_i(t)=0\;(i=1,\ldots,N).$$  
\end{proposicao}

It is worth to remark that  $S(t_0)<{S^*}$ and ${Z}={S^*}$ in equation (\ref{Uliapunov}) leads to ${\dot U}<0$ for all $t>t_0$ in equation (\ref{dUdt1}). Then, $U(t)$ is a strictly decreasing function for all $t>t_0$ and we have 
\begin{equation}
\label{woc}
U(t)< U(t_0),\;\;\forall t>t_0
\end{equation}
The inequality above is called {\it weak outbreak condition}.
Since $S(t)$ is a non-increasing function, we can establish the following result:

\begin{proposicao}(The Weak Outbreak Condition) The weak outbreak condition is met for every $ t> t_0 $ such that $ S_0 = S (t_0) <{S^*} $.\end{proposicao}

Finally,  from the expressions for ${\bar S}$ and ${S^*}$, respectively  in (\ref{SO}) and (\ref{SL}), it is straightforward to show that ${\bar S}<{S^*}$.

\section{Mathematical Definition of a Control Strategy}

The system of differential equations in (\ref {seir}) is an autonomous system as long as the parameters that appear in the differential equations are independent of time. In order to define an epidemic control strategy, we may consider this system of equations as non-autonomous provided the basic reproduction numbers $ R_i (t) \; (i = 1, \ldots, N) $ are modeled by time-dependent functions.
We name the vector function $ {\vec R} (t) = (R_1 (t), \ldots, R_N (t)) $ as the basic reproduction vector, where the functions $ R_i (t) $ are defined in the range $ 0 \leq t <\infty $. The initial value $ {\vec R}^0 $ (or natural value) of the basic reproduction vector is its value at time $ t = 0 $, that is, $ {\vec R}^0 = (R_1 (0 ), \ldots, R_N (0)) $. 

A control strategy is modeled by a basic reproduction vector that is given generically as
\begin{equation}
\label{controle}
{\vec R}(t)=\left\{
\begin{array}{ll}
{\vec R}^0 & 0\leq t <t_I\;\; \\  & \\
{\vec R}^c(t) &  t_I\leq t<\infty
\end{array}\right.,
\end{equation}
where $ {\vec R}^c (t) $ is a well-defined vector function for every $ t \geq t_I $ with $ t_I> 0 $. The instant $ t_I $ is the initial moment when the control strategy begins and the vector function $ {\vec R}^c (t) $ is called the controlling function of the basic reproduction numbers. Therefore, we must consider the non-autonomous system (\ref{seir}) taking into account that  $\beta_i(t)$ are time-dependent functions given by
\begin{equation}
\beta_i(t)=\gamma_i R_i(t)\;(i=1,\ldots,N).
\end{equation}

A convenient way to characterize a control strategy is to define the following vector:
\begin{equation}
{\vec q}(t)=\left(q_1(t),\ldots,q_N(t)\right)=\left(\frac{R_1(t)}{R^0_1},\ldots,\frac{R_N(t)}{R^0_N}\right).
\end{equation}
In this way, any control strategy given in (\ref{controle}) can also be defined as
\begin{equation}
q_i(t)=\left\{\begin{array}{ll}
1 & 0\leq t \leq t_I \\ & \\
\displaystyle \frac{R_i^c(t)}{R_i^0} & t_I< t <\infty
\end{array}\right.,\;\;(i=1,\ldots,N).
\end{equation}
Also, the values of ${\bar S}$ and ${S^*}$, respectively defined in (\ref{SO}) and (\ref{SL}), become time-dependent functions written as  
\begin{equation}
\label{SLcontrol}
{\bar S}(t)=\left(\displaystyle\sum_{i=1}^N q_i(t)\frac{x_i\gamma_i}{\gamma} R^0_i\right)^{-1},\;\;
{S^*}(t)=\left(\displaystyle\sum_{i=1}^N q_i(t) x_i R^0_i\right)^{-1}.
\end{equation}
The values of the functions above evaluated at $t=0$  are associated to the epidemic evolution without control and they are obtained plugging $q_i(t)=1$ in (\ref{SLcontrol}), which leads to
\begin{equation}
\label{SLnatural}
{\bar S}_0=\left(\displaystyle\sum_{i=1}^N \frac{x_i\gamma_i}{\gamma} R^0_i\right)^{-1},\;\;
S^*_0=\left(\displaystyle\sum_{i=1}^N x_i R^0_i\right)^{-1}
\end{equation}
 
A strategy is called permanent if $ \lim_{t \rightarrow \infty} {\vec R}^c (t) \neq {\vec R}^0 $, that is, the basic reproduction vector never returns asymptotically to its natural (or initial) value or, equivalently, for some component $i$ of the vector ${\vec q}(t)$ we have $\lim_{t\rightarrow\infty}q_i(t)\neq 1$ . 
A strategy is said to be non-permanent if  $ \lim_{t \rightarrow \infty} {\vec R}^c (t) = {\vec R}^0 $ or, equivalently, $\lim_{t\rightarrow\infty}q_i(t)= 1$ for all $i=1,\ldots,N$. If there is an instant of time $ t_F> t_I $ such that $ {\vec R}^c (t) = {\vec R}^0 $ for all $ t> t_F $, then the strategy is called finite  and its time duration is defined as $ \Delta t = t_F-t_I $.

Let us now characterize two important kind of strategies:
\begin{itemize}
\item{\it Non-Outbreak Strategies} (NOS), where $S(t)\leq {\bar S}(t)$ for all $t>t_I$;
\item {\it Weak Outbreak Strategies} (WOS), where  $S(t)\leq S^*(t)$ for all $t>t_I$.
\end{itemize} 
It is easy to see that for non-outbreak strategies, the solution of the respective  non autonomous system (\ref{seir}) satisfies the non-outbreak condition established in proposition 1 
for all $t>t_I$.

On the other hand, let us plug ${Z}=S^*(t)$ and $\beta_i=\gamma_i R_i(t)$ in (\ref{Uliapunov}) to obtain:
\begin{equation}
\label{Ut}
U(t)=E(t)+\sum_{i=1}^N{S^*(t)R_i(t)}I_i(t)=E(t)+\sum_{i=1}^N\left(\frac{R_i(t)}{\sum_{j=1}^Nx_jR_j(t)}\right)I_i(t)
\end{equation} 
and 
\begin{equation}
\label{dUt}
{\dot U}(t)=\left(S(t)-S^*(t)\right)\sum_{i=1}^N\gamma_i R_i(t)I_i(t)+\sum_{i=1}^N\frac{d}{dt}\left(\frac{R_i(t)}{\sum_{j=1}^Nx_iR_i(t)}\right)I_i(t).
\end{equation}

A strategy is said to be uniform if $ q_i (t) = q (t) $ for all $ i = 1, \ldots, N $, that is, in a uniform strategy all the components of the basic reproduction vector are always changed in the same proportion and the proportionality constant is given by a certain function $ q (t) $. If we substitute $R_i(t)=q(t)R_i^0$ in equations (\ref{Ut}) and (\ref{dUt}) we obtain
\begin{equation}
\begin{array}{c}
\displaystyle U(t)=E(t)+\sum_{i=1}^N\left(\frac{R^0_i}{\sum_{j=1}^N x_j R_j^0}\right)I_i(t) \\ \\
\displaystyle {\dot U}(t)=\left(S(t)-S^*(t)\right)q(t)\sum_{i=1}^N\gamma_i R^0_i I_i(t).
\end{array}
\end{equation} 
Hence, for any uniform WOS  we have ${\dot U}\leq 0$  for all $t>t_I$ and, consequently, $U(t')\leq U(t)$ for all $t'>t>t_I$. This property justifies the name Weak Outbreak for this kind of strategy.  

The main point is to define and analyze the mitigation effects in the spread of an epidemic, depending on the application of a given control strategy that alters the natural course of its development.
Regarding the analysis of the effects of a finite non-permanent control strategy, we have to keep in mind their two main objectives: (i) to reduce the total number of infected individuals after the end of the epidemic cycle and (ii) decrease the peak intensity of the epidemic. Naturally, these two objectives must be evaluated according to the duration of the control strategies, since extremely long duration can entail very high economic and social costs.

Regarding objective (i), we should analyze the asymptotic value of susceptible $ S_\infty = \lim_ {t \rightarrow \infty} S (t) $, since the total number of infected individuals after the end of the epidemic cycle is given by $ 1-S_\infty $.
From proposition 1, it is straightforward to proof an extremely important result for finite strategies that are applied in models described by the system of differential equations given in (\ref{seir}):

\begin{proposicao} (Fundamental Principle of Finite  Strategies)
For any finite strategy applied in a model given by (\ref{seir}) we have
$$S_{\infty}=\lim_{t\rightarrow\infty}S(t)< S^*_0,$$
where $ {S^*_0} $ is defined in the equation (\ref{SLnatural}).
\end{proposicao}

The main consequence of the fundamental principle is that for any finite strategy the total number of infected individuals is always greater than $ 1-{S^*} $. This result implies a limitation for finite  strategies in order to diminish the final number of infected. This happens because it is not possible to control the value of $ {S^*_0} $, associated with the natural reproduction vector $ {\vec R}^0 $ of the epidemic.
On the other hand, permanent  strategies can increase the value of $ {S^*_0} $ and, consequently, to lower the minimum limit
for the total number of infected, however, there is a price to pay: to ensure a permanent control.

These properties appear very clearly in the applications that one  makes using uniform strategies in a SEIAR model and non-uniform strategies, associated to the isolation of  infected individuals, in a SEIAQR model.

\subsection{Application to a Model where the $N$ Parallel Stages have the same infection time} 

Let us consider the model described in (\ref{seir}) with the additional condition that the infection times are all the same, that is,
$ \gamma_i = \gamma \; (i = 1, \ldots, N) $. Without loss of generality we write $ R_i = e_iR_0 $, where $ R_0 $ is an arbitrary reproduction rate used as a reference value, and we define $ \beta_i =e_i \beta $ with $ \beta = \gamma R_0 $ . In this way, the model (\ref{seir}) can be rewritten as:
\begin{equation}
\label{seiarap}
\begin{array}{c}
{\dot S}=-\beta\left(\sum_{i=1}^{N}e_iI_i\right)S,\;
{\dot E}=\beta\left(\sum_{i=1}^{N}e_iI_i\right)S-\sigma E, \\ \\
{\dot I}_i=x_i\sigma E-\gamma I_i,\;
{\dot R}=\gamma\left(\sum_{i=1}^NI_i,\right),
\end{array}\;(i=1,\ldots,N)
\end{equation}
where $ e_i $ is the infectiousness of  infected  in stage $ i $ relative to the reference value $ R_0 $.
We can also conclude that $ {\bar S} $ and $S^*$, respectively given in (\ref{SO}) and (\ref{SL}), have the same value

\begin{equation}
\label{SLap}
{\bar S}={S^*}=\left(\sum_{i=1}^N x_i e_i\right)^{-1}R_0^{-1},
\end{equation} 
and then any WOS is a NOS.   

Considering the variable $J=\sum_{i=1}^N e_iI_i$ we can show that $S$, $E$ and $J$ obey the following three dimensional system of differential equations:
\begin{equation}
\label{seired}
{\dot S}=-\beta JS,\; 
{\dot E}=\beta JS-\sigma E,\;
{\dot J}={\bar\sigma} E-\gamma J\;\left({\bar\sigma}=\sigma\sum_{i=1}^Nx_ie_i\right), 
\end{equation}
which is decoupled from the equations for $I_i\;(i=1,\ldots,N)$ and $R$  in the system  (\ref{seiarap}).
We can now state a condition for which the value of $ {S^*} $ in (\ref{SLap}) may be interpreted as a critical value for solutions of system (\ref{seired}). This property, which is shown in appendix A, can be proved for solutions with initial condition that implies an outbreak of the epidemic.
We say that an initial condition in an instant $ t_0 $ (given by $ S_0 = S (t_0) $, $ E_0 = E (t_0) $ and $ J_0 = J (t_0) $) is an initial condition of outbreak if there is a time $ t_1> t_0 $ for which $ {\dot E} (t_1)> 0 $ and $ {\dot J} (t_1)> 0 $. From this definition we can proof the proposition (see appendix A):

\begin{proposicao}
(The Critical Value of Epidemic) For any solution of the differential system (\ref{seired}) with initial condition of outbreak, there is a time $ t^*> t_1 $ such that
\begin{eqnarray}
&&{\dot E}(t)>0,\;{\dot J}(t)>0\;{\rm for}\; t_1<t^*,\nonumber \\
&&{\dot E}(t)=0,\;{\dot J}(t)=0\;{\rm for}\; t=t^*, \nonumber\\
&&{\dot E}(t)<0,\;{\dot J}(t)<0\;{\rm for}\; t>t^* , \nonumber
\end{eqnarray}
where $S(t^*)={S^*}$.
\end{proposicao}

\section{SEIAR Model for COVID-19: an Analysis of Uniform Control Strategies }

The SEIAR model for COVID-19 is obtained from the system of ODE's  in (\ref{seiarap}) assuming the existence of two stages of infected: the symptomatic infected $ I_1 = I $ and the asymptomatic infected $ I_2 = A $. The basic reproduction number for $ I_1 $ and $ I_2 $ are respectively defined as $ R_1 = R_0 $ and $ R_2 = \xi R_0 $, where $ \xi $ is the infectiousness of asymptomatic relative to symptomatic. Therefore, the corresponding system of differential equations in (\ref{seiarap}) for $N=2$ becomes
\begin{equation}
\label{seiar}
\begin{array}{ll}
{\dot S}=-\beta (I+\xi A) S,\;
{\dot E}=\beta (I+\xi A) S-\sigma E,\\ \\
{\dot I}=(1-\chi)\sigma E-\gamma I,\;
{\dot A}=\chi\sigma E-\gamma A,\;
{\dot R}=\gamma (I+A),\;
\end{array}\;\beta=\gamma R_0.
\end{equation}

We summarize the description of variables, parameters and their respective empirical values for COVID-19 in table 2 below. The empirical values  were based on the values obtained in references \cite{ref14}, \cite{ref15}, \cite{ref16}, \cite{ref18}, \cite{ref39} and  are chosen in order to define the ENM associated to the city of Brasilia-Brazil, which has a resident population of $ 3 $ million.

\vspace*{0mm}
\begin{center}
\begin{tabular}{c}
Table 2 - Variables and Parameters in ODE's system (\ref{seiar})\\
\hline\hline 
Variables of the SEIAR Model \\
\hline\hline 
\begin{tabular}{cc}
$S$ & Susceptible Population\\ \\
$E$ & Exposed Population \\ \\
$I$ & Symptomatic Infected Population \\ \\
$A$ &  Asymptomatic Infected Population\\ \\
$R$ &  Recovered Population
\end{tabular} \\
\hline\hline
Parameters of the SEIAR Model for COVID-19 \\ \hline\hline 
\begin{tabular}{ccc}
$R_0$ & Basic Reproduction Rate of Symptomatic Infected &  $3.0$ \cite{ref38} \\  \\
$\sigma^{-1}$ & Incubation Time & 5.0 days \cite{ref14}\\ \\
$\gamma^{-1}$ & Infection Time & 1.61 days \cite{ref14}\\ \\
$\chi$ & Likelihood of an Exposed becomes Asymptomatic & 0.862 \cite{ref15} \\ \\
$\xi$ & Asymptomatic Infectiousness & 0.55 \cite{ref15}\\
&Relative to Symptomatic & \\
\hline\hline
\end{tabular}
\end{tabular}
\end{center}

Our objective in this section is to study the main properties of some uniform strategies applied to a  SEIAR model and  illustrate them trough the numerical integration of system (\ref{seiar}), using the empirical values for the parameters. 
It is worth to remark that we are not interested in to assess whether the chosen ENM, with empirical  parameters given in table 2, are the best to describe the evolution of the epidemic in the city of Brasília, but rather to analyze the properties of some control strategies.
We seek to make our analysis as general as possible and not strictly dependent on the specific parameters used in the numerical integration.

\begin{figure}
[!htb]
\begin{center}
\includegraphics[width= 14.0cm]{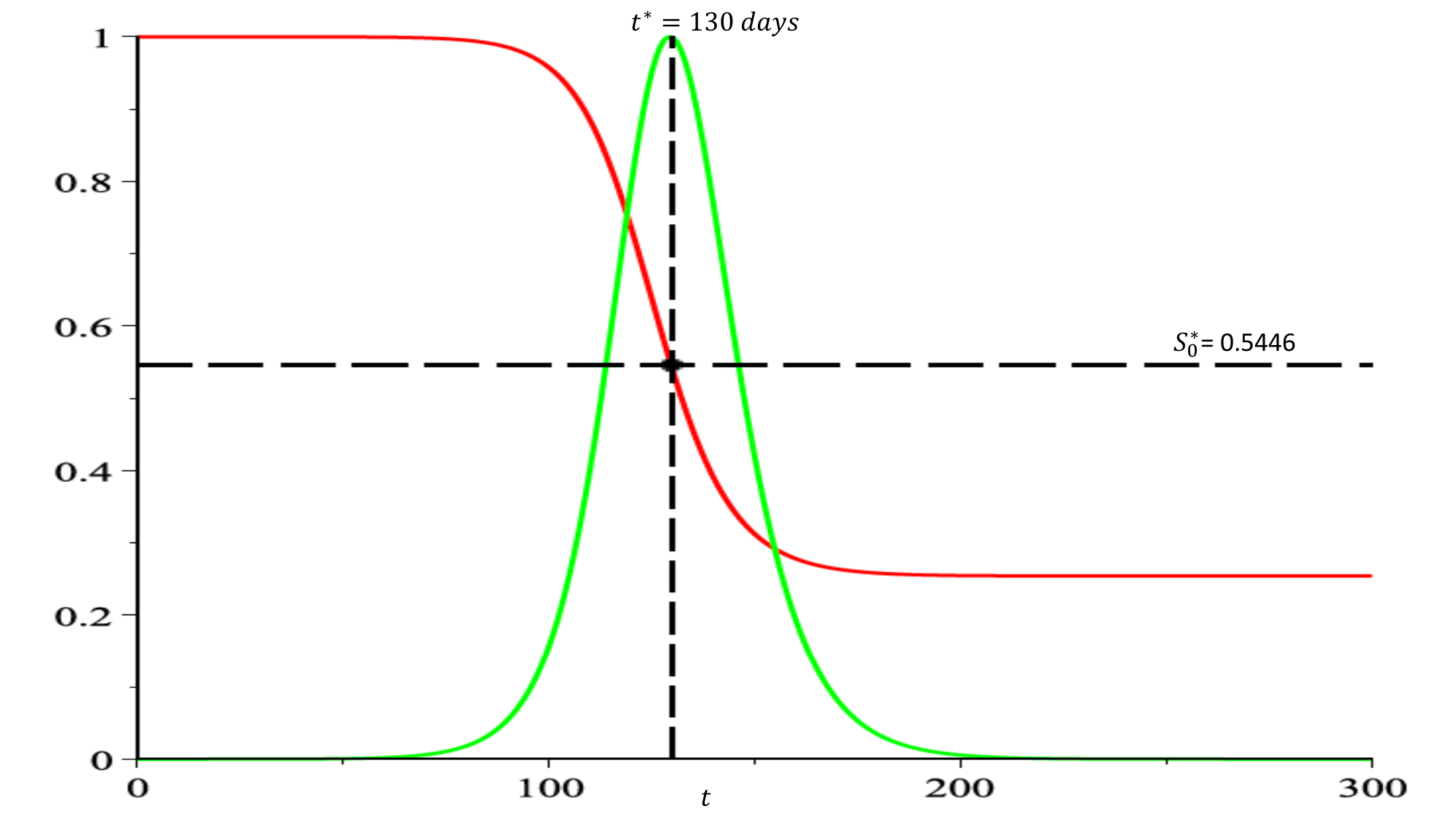}
\end{center}
\vspace*{-9mm}
\caption{ $ E (t) / E_ {max} $ (Green Curve), $ S (t) $ (Red Curve). The critical value $ {S_0^*}$ correspond to the horizontal dashed lines. The vertical dashed line represents the time $t^*$ such that $S(t^*)={S_0^*}$. $ E_ {max} $ is the maximum value of the function $E(t) $.}
\end{figure}

Before to study uniform $R_0$-Strategies, we present  the solution of the ENM for Brasilia given in (\ref{seiar}) with parameter values in the table 2.
We consider an initial condition with a single exposed $ E (0) = 1 / N $, where $N=3\times10^6$.
Figure 1 shows the time evolution of  variables $ S (t) $ and $  E (t)  $.
We can clearly see how the  value $ {S_0^*} $  represents the inflection in the epidemic evolution and corresponds to the maximum value of  $ E (t) $. The critical value is reached at approximately $ t^* = 130 $ days  when $ {S_0^*} = S (t^*) $. This figure illustrated very well the content of proposition 5.
Indeed, the curves presented in figure 1 represents the natural evolution of  the epidemic, that is, the epidemic evolution in the absence of any control strategy.

The  corresponding natural value $ {S^*_0} $ for the SEIAR model in (\ref{seiar}) is obtained from  equation (\ref{SLap}), substituting the parameter values in table 2:
\begin{equation}
\label{SLSH}
\displaystyle {S^*_0}=\frac{1}{1-\chi+\xi\chi}\frac{1}{R_0}=0.5445,
\end{equation}

\subsection{Uniform $R_0$-Strategies}

An uniform control strategy, denoted as $R_0$-Strategy, essentially means to redefine the value of $ R_0 $, hence the value
of $ \beta $ in the SEIAR system given  in (\ref{seiar}). Therefore,
a $R_0$-Strategy  is mathematically defined  as follows:
\begin{equation}
\label{estrategia}
\beta=\gamma R_q,\;\;\; R_q=q(t)R_0,\;\;\;
\left\{\begin{array}{ll}
q(t)=1 & 0<t\leq t_I \\ \\
q(t)\leq 1 & t>t_I
\end{array}\right.,
\end{equation}
where there is a $\delta>0$ such that $q(t)<1$ for all $t_I<t<t+\delta$. The time $ t_I $ is the initial moment in which the control strategy starts to be applied and the SEIAR system (\ref{seiar}), with the new definition of $ \beta $, becomes a non-autonomous differential system.

In this section we define some basic concepts associated with  $R_0$-Strategies and illustrate them through the SEIAR system in (\ref{seir}) used to model the control strategy applied to the city of Brasilia. In all numerical integration, we consider the empirical parameters in table 2  and an initial condition with a single exposed $ E (0) = 1 / N $, where $ N = 3 \times 10^6$.

First, we define the critical susceptible function $S_q(t)$, associated to a given $R_0$-Strategy, as being the function ${S^*}(t)$  given in (\ref{SLcontrol}). Hence, we have
\begin{equation}
\label{Sq}
S_q(t)={S^*}(t)=\frac{1}{(1-\chi+\xi\chi)}\frac{1}{R_q(t)}=\frac{1}{q(t)}\frac{1}{(1-\chi+\xi\chi)}\frac{1}{R_0}=\frac{{S^*_0}}{q(t)}
\end{equation} 
and, by the definition of function $ q (t) $,
\begin{equation} 
S_q(t)= {S^*_0}\;\;\;  (t\leq t_I);\;\;\;\;
S_q(t)\geq {S^*_0}\;\;\;  (t>t_I).
\end{equation} 
Next, we retake the classification of some types of control strategies given in section 2 in the particular context of $R_0$-Strategies applied to a SEIAR model.
\begin{itemize}
\item {\it Permanent Strategy} is a strategy where $ \lim_ {t \rightarrow \infty} q (t) <1 $.
\item {\it Non-Permanent Strategy} is a strategy where $ \lim_ {t \rightarrow \infty} q (t) = 1 $.
\item {\it Finite Strategy} is one in which there is an instant of time $ t_F> t_I $ such that $ q (t) = 1 $ for all $ t> t_F $.
\end{itemize}

An important point to analyze in finite strategies is their duration: $ \Delta t = t_F-t_I $. The shorter this duration, the lower the burden of reducing the basic reproduction number, implying less economic costs. Let us remark that finite strategies have an important  practical interest, since the existence of an end time  for the control strategy implies that the normal situation returns after this time. On the other hand, the permanent strategies are associated with permanent change in the the basic reproduction number of the epidemic.

We can clearly see that a $R_0$-Strategy such that $ q (t) \leq {S_0^*} / S (t) $ for all $ t> t_I $ is a non-outbreak strategy. Indeed, in  such kind of  $R_0$-strategy, the redefined value of $R_0$ implies that the dynamical system SEIAR is always in a non-outbreak condition. 

\subsection{Constant Strategies}

The simplest strategies are those in which the value of $ q (t) $  changes at most twice and we call them as {\it Constant Strategies}. There are two types of constant strategies:
\begin{itemize}
\item {\it Constant Permanent Strategy} where the $ q (t) $ function is defined as
\begin{equation}
q(t)=\left\{\begin{array}{ll}
1 & 0<t\leq t_I \\ \\
q_I & t>t_I
\end{array}\right.\; ;
\end{equation}
\item {\it Constant Finite Strategy} where the $ q (t) $ function is defined as
\begin{equation}
q(t)=\left\{\begin{array}{ll}
1 & 0<t\leq t_I \\ \\
q_I & t_I<t<t_F \\ \\
1 & t\geq t_F
\end{array}\right.\; ;
\end{equation}
\end{itemize}
where $ q_I $ is a constant value such that $ q_I <1 $. Therefore,  a {\it Constant Permanent Strategy} is a non-outbreak strategy if  
$ q_I \leq {S_0^*} / S (t_I) $ or $ q_I \leq {S_0^*} $.

\begin{figure}
[!htb]
\begin{center}
\includegraphics[width= 14.0cm]{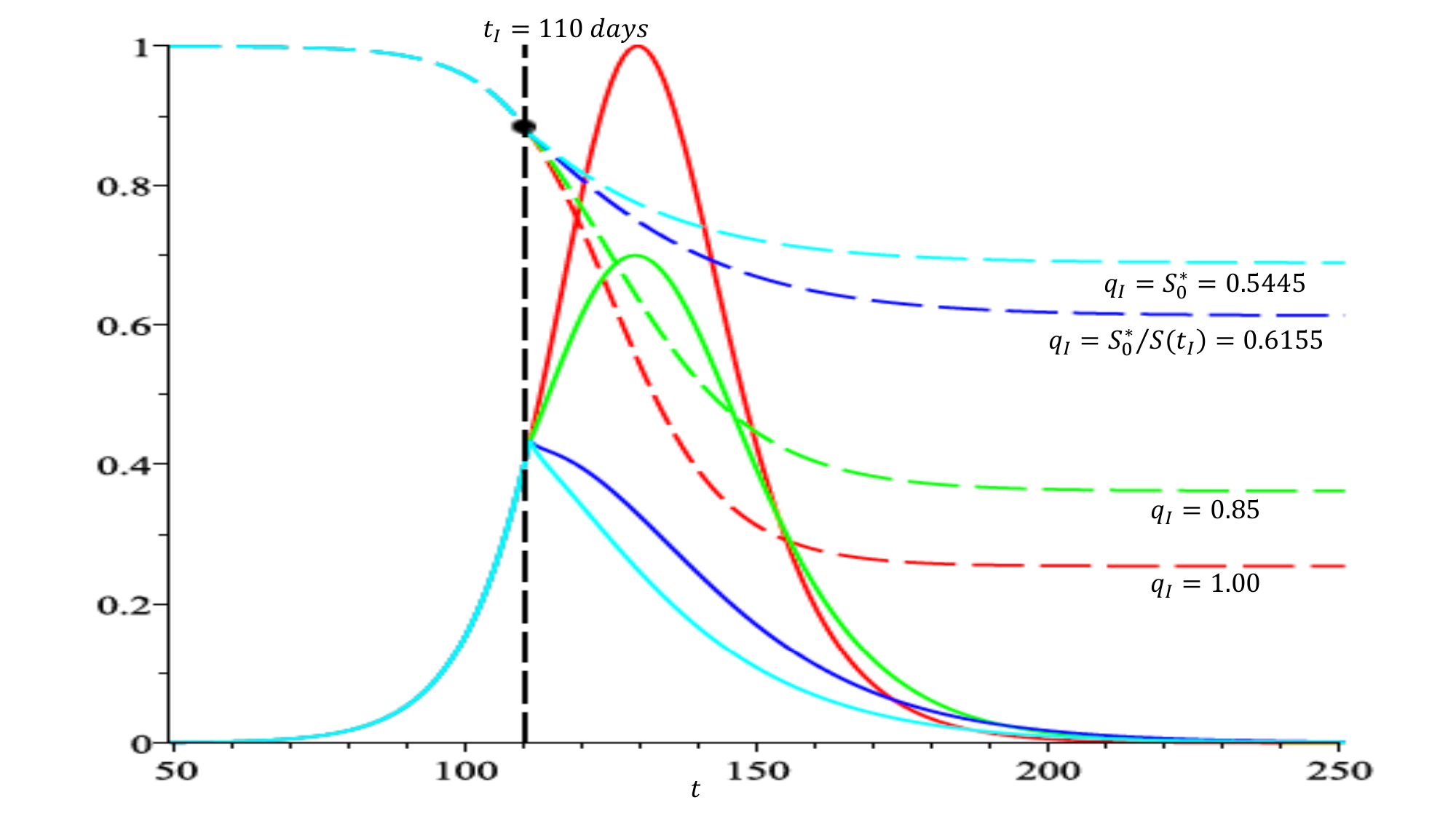}
\end{center}
\vspace*{-8mm}
\caption{$||{\vec V(t)}|| / V_ {max} $ (Solid Curves) and $ S (t) $ (Dashed Curves).  The vertical dashed line represents the start time $ t_I = 110 $. $ ||{\vec V}(t)||_ {max} $ is the maximum value of  $ ||{\vec V}(t)|| $  for $ q_I = 1.0 $. Each color corresponds to a strategy and the respective values of $ q_I $ are shown in the figure.}
\end{figure}

Figure 2  illustrates the time evolution of three {\it Constant Permanent Strategies}, defined by different values of $ q_I $. The evolution of the system submitted to these strategies is compared with the natural epidemic evolution ($ q_I = 1 $) without control strategy.
We can see, respectively for the values  $ q_I = {S_0^*} / S (t_I) $ and $ q_I = {S_0^*} $, that the corresponding $R_0$-Strategies are non-outbreak, whereas the strategy with $ q_I = 0.85> {S_0^*} / S (t_I) $ is an outbreak strategy.

Figure 3 shows the time evolution of four {\it Constant Finite Strategies}, all defined by the same value of $ q_I = {S_0^*} $, starting  in $ t_I = 110 $ days and lasting different time intervals ($ t_F = t_I + \Delta t $ with $ \Delta t = 30, 60, 90, 120 $ days). The evolution of the system submitted to these strategies is compared with the natural evolution of the epidemic ($q_I=1$)
and with the permanent strategy starting at $ t_I$  ($ q_I = {S_0^*} $).
We note  that for all finite strategies the asymptotic value of susceptible is less than $ {S_0^*} $, that is, $ \lim_ {t \rightarrow \infty} S (t) <{S_0^*} $, 
which illustrates the content o proposition 4.
Also, we emphasize that all finite strategies shown in Figure 3 are outbreak strategies for $t>t_F$, this is clearly due to the fact that they do not satisfy the condition $ S (t_F)  \leq {S_0^*} $. 
Let us observe that while the finite strategy is applied it follows the trajectory of the permanent strategy defined with the same value of $ q_I $, so we can never have $ S (t_F) <{S_0^*} $ for any final time $t_F>t_I$.

\begin{figure}
[!htb]
\begin{center}
\hspace*{-2mm}\includegraphics[width= 14.0cm]{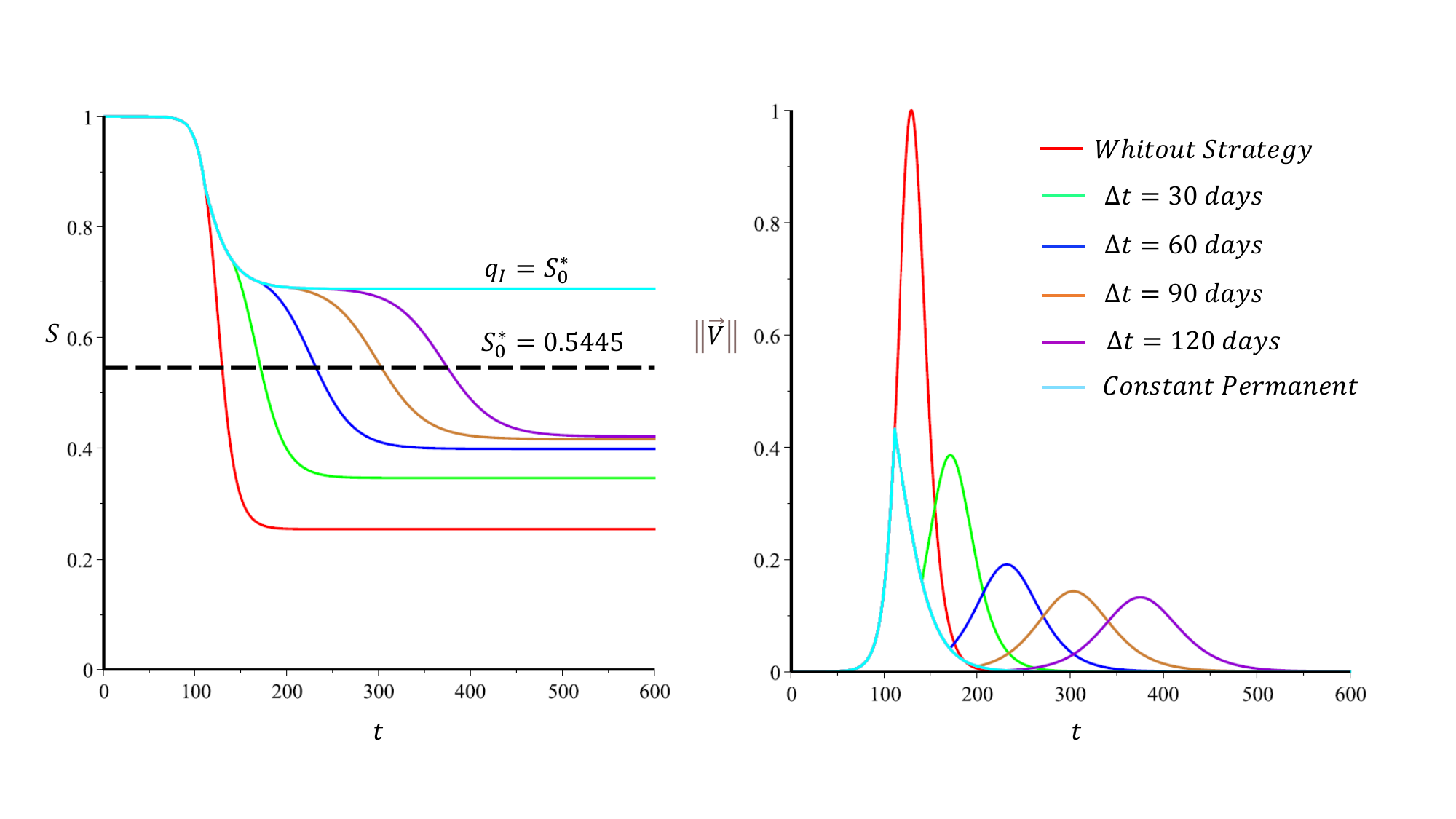}
\end{center}
\vspace*{-10mm}
\caption{Comparison of finite constant strategies through the time evolution of $  E (t)  $ (Right Panel) and $ S (t) $ (Left Panel). All strategies start at $ t_I = 110 $ days and are defined by the same value as $ q_I $. Finite strategies have a duration of $ \Delta t = t_F-t_I $. }
\end{figure}

All constant finite strategies showed in figure 3 have a time $t_F$ where the value of $R_0$ returns to its initial value and they display a similar pattern: they have a clear second peak that tends to have the same shape as times $t_F$ increases (see the right panel in figure 3).
Indeed, this indicate the irrelevance to keep the control strategy after some time $\Delta t_c$ since, whatever the time $t_F>t_I+\Delta t_c$, 
the second peak will have approximately the same shape.
More specifically, we can estimate $\Delta t_c\approx 90$ days, so that for a time $t_F>t_I+90$ days (see the gold and violet curves) the asymptotic value of $S(t)$ and the peak of $||{\vec V}(t)||$ are approximately the same. Let us observe (in the left panel) how the asymptotic value of $S(t)$ tends to be the same for the gold  ($\Delta t=90$) and violet ($\Delta t=120$) curves.

\begin{figure}
[!htb]
\begin{center}
\hspace*{-2mm}\includegraphics[width= 14.0cm]{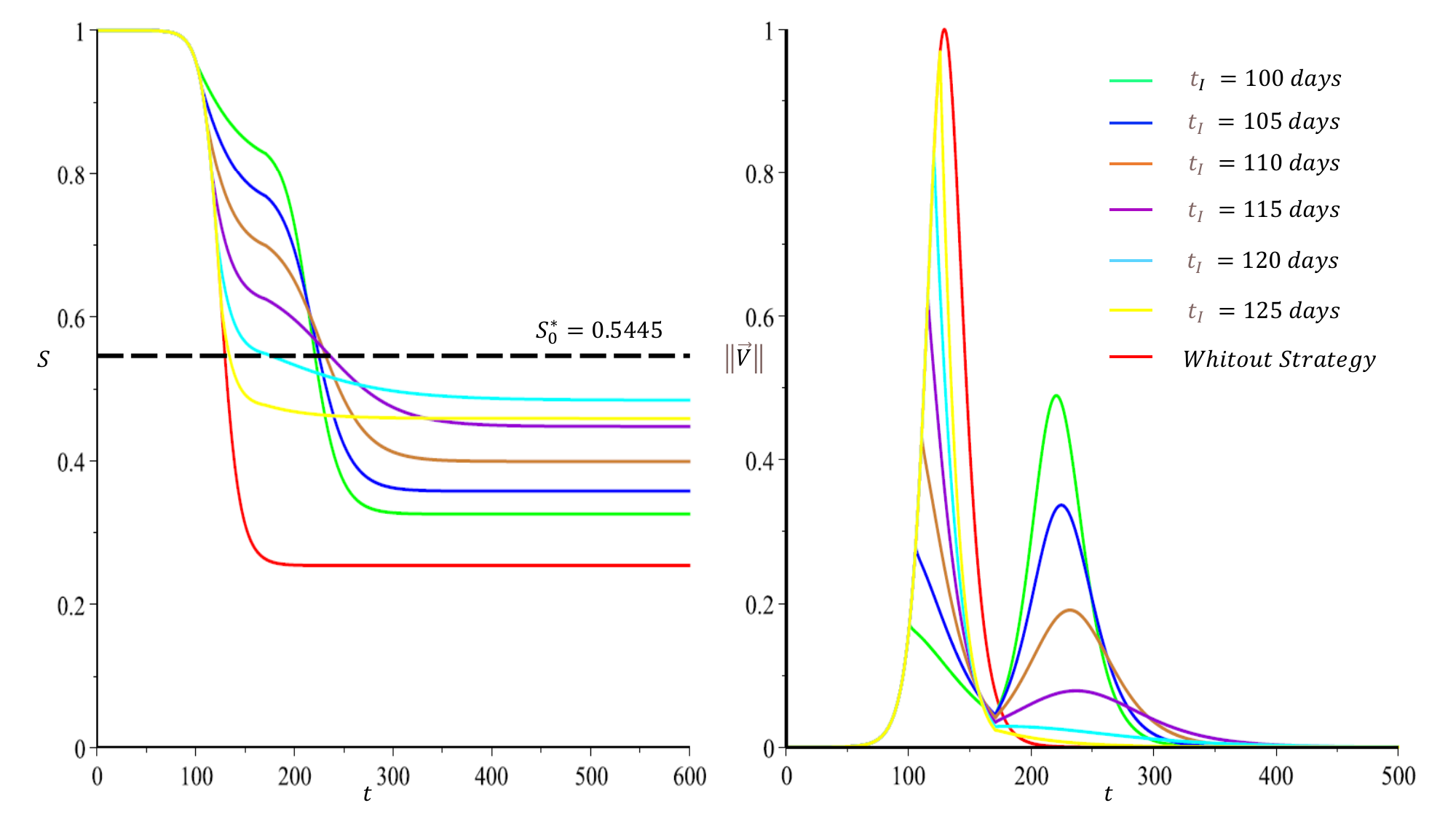}
\end{center}
\vspace*{-10mm}
\caption{Comparison of finite strategies through the time evolution of $ || \vec V (t) || $ (Right Panels) and $ S (t) $ (Left Panels). All strategies finish at $ t_F = 170 $ days and are defined by the same value of $ q_I=S_0^* $. The strategies start at different times  $ t_I $ corresponding to graphs of different colors. }
\end{figure}

\begin{figure}
[!htb]
\begin{center}
\hspace*{-2mm}\includegraphics[width= 14.0cm]{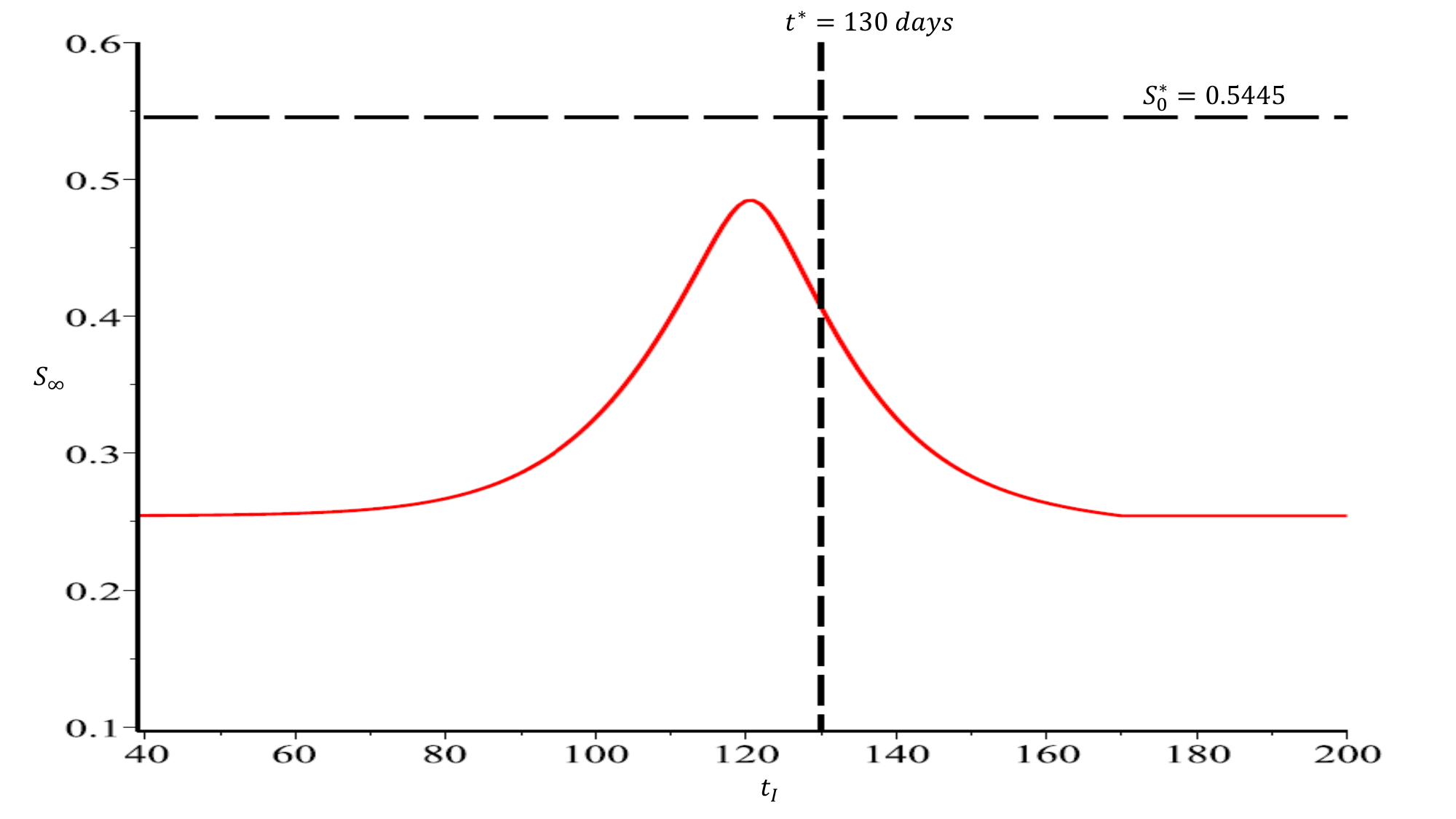}
\end{center}
\vspace*{-10mm}
\caption{Asymptotic value  $S_{\infty}$ as a function of $t_I$ for Finite Strategies }
\end{figure}

In figure 4 we compare different constant strategies that end in the same final time $ t_F> t^* $, but start at different initial times $ t_I <t^* $ (remember that $ t^*$ is the time in which the number of incubates $ E (t) $ reaches its peak in the natural evolution of the epidemic).
In general, we can conclude that all strategies break the single peak of the epidemic in its natural evolution into two peaks separated by a certain period of time.
Also, it is quite clear that the earlier the strategy starts, the greater the reduction in the first peak and the higher the second peak in relation to this one.

An important question is to assess the asymptotic value of susceptible $ S_\infty = \lim_ {t \rightarrow \infty} S (t) $ that allows us to calculate the total number of infected after the end of the epidemic cycle. In figure 5, we show the asymptotic value of susceptible $S_\infty$  as a function of time $ t_I $.
We can see that finite constant strategies have some significant effect on the asymptotic value of susceptible only in the range $ 90 <t_I <150 $.
In other words, we see that starting to apply the control strategies too soon, or too late, does not mean reducing significantly the 
total number of infected individuals after the end of the epidemic cycle. It is worth mentioning that constant strategies have a significant effect on $ S_\infty $ even when they start after the peak of the epidemic ($ t_I> t^* $).
From figure 5 it is clear that there is an optimal initial time $ t_I $ to start applying the control strategy. This time (corresponding to the maximum point of the red curve showed in the figure) is approximately $ t_I = 120 $, that is, the optimum time to start the control strategies is approximately 10 days before the peak ($ t^* = 130 $).

It is important to draw attention to the fact that the optimal time was defined as the one in which the final number of susceptible is the greatest, implying the smallest number of infected  at the end of the epidemic cycle. However, this is not the only important criterion to decide whether a strategy is better or worse. We must remember that it is very important to take into account the maximum value reached for the epidemic, as an high peak value can mean that many individuals are infected at the same time and, consequently, the healthcare system would be stressed and not be able to treat all those infected individuals with severe symptoms, which can lead to a higher mortality rate.

In this way, we clearly see a trade-off problem between these two ways of assessing the effectiveness of a strategy: the smallest final number of infected versus the smallest number of infected at the peak of the epidemic. With regard to what was shown in figures 4 and 5, we see that the optimal time to decrease the final number of infected $ 1-S_\infty $ corresponds to an intervention that should start very close to the peak of the epidemic, implying in starting control too late, in a situation where the number of infected in a short period of time is very high.

The biggest problem with the constant finite strategies is that when they finish the non-outbreak condition may not be satisfied, implying a second outbreak of the epidemic.
An important problem here is that  stronger and earlier the intervention happens (trying to avoid the natural peak of the epidemic), the greater the chance that the non-outbreak condition at the end of the strategy is not fulfilled.
In this sense, it would be interesting to define finite control strategies that always satisfy the non-outbreak condition.

\subsection{Weak Outbreak Strategies (WOS)}

A regulated strategy is a control strategy where the value of $q(t)$ depends on the knowledge of some information  about the actual state of the SEIAR model in (\ref{seir}). We now define the following  regulated strategy:
\begin{itemize}
\item $S_0^*$-Strategy is a strategy that has a function $ q (t) $ given by
\begin{equation}
q(t)=\left\{\begin{array}{cl}
1 & 0<t\leq t_I \\ \\
\left\{\begin{array}{cl}
\displaystyle \frac{S_0^*}{S(t_i)} &  S(t_i)>S_0^* \\
1 &  S(t_i)\leq S_0^*
\end{array}\right.
& t_i<t<t_{i+1}
\end{array}\right.
\begin{array}{c}
t_0=t_I,\;t_{i+1}-t_i=\Delta\\
(\Delta>0,\;\;i=0,\ldots\infty)
\end{array}.
\end{equation}
\end{itemize}
Let us remark that in this kind of strategy the value of $ q (t) $ implies that $ S_q (t) \geq S (t) $ for all $t$ (see definition of $ S_q (t) $ in the equation (\ref{Sq})) and we can conclude that 
a $S_0^*$-Strategy is a WOS. Also, for a model given in (\ref{seiarap}) follows the relation (\ref{SLap}), which draw in to the conclusion that  a $S_0^*$-strategy fulfills the non-outbreak condition for all $t>t_I$. Indeed, for a model described by the system in (\ref{seiarap}) a WOS is equivalent to a NOS.

The parameter $\Delta>0$ is called period and
we illustrate in figure 6 two $S^*_0$-Strategies: one with period $\Delta=1$ day and the other with $\Delta=30$ days.
Both strategies start  at $ t_I = 110 $ days and, for comparison, we plot the curves of natural epidemic evolution ($q(t)=1$) and with constant permanent strategy ($ q_I = S_0^* $).

\begin{figure}
[!htb]
\begin{center}
\includegraphics[width= 14.0cm]{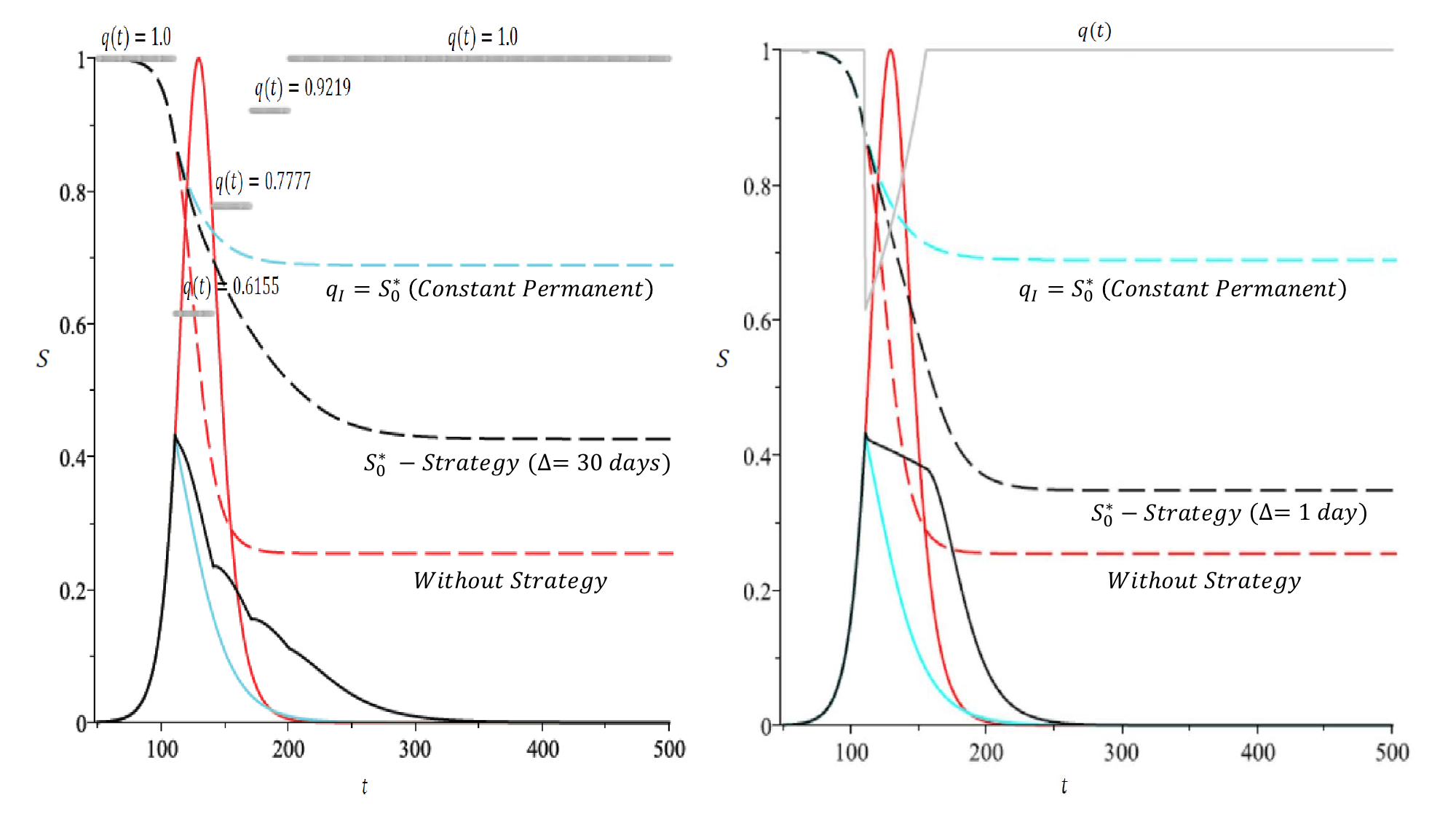}
\end{center}
\vspace*{-10mm}
\caption{ $ ||{\vec V} (t) || / V_ {max} $ (solid curves) and $ S (t) $ (dashed curves).
$ ||{\vec V}(t)||_ {max} $ is the maximum value of $  ||{\vec V}(t)||  $ for the epidemic without strategy. The $ q(t) $ function  is the gray curve. Each color represents a different strategy. All the strategies start at $ t_I = 110 $ days.}
\end{figure}

The first important remark is that $S_0^*$-Strategies are finite, implying a well-defined end time $ t_F $. The second important remark relates to the fact that a $S_0^*$-Strategy prevents the existence of a second peak after the end of the control. This stems from the fact that the strategy, by the way it is defined, only ends when the proportion of susceptible falls below the critical value $ {S_0^*} $ and, from that moment on, the epidemic is already in a natural condition of non-outbreak.
A $S_0^*$-strategy clearly flattens the peak of the epidemic and implies a lower number of asymptotic infected ($1-S_\infty$) than in the case of the natural evolution of the epidemic without control. So, it is interesting to investigate, in the same way as done for finite constant strategies, how the final (asymptotic) number of susceptible depends on the start time $ t_I $.

In figure 7 we show the effect of starting a $S_0^*$-strategy at different times $ t_I $.
From the left panel of this figure, we see that the earlier the control strategy begins, the longer it lasts and the longer and flatter the peak of the epidemic. We also see (in the smaller left panel) that the sooner control begins, the greater the asymptotic of susceptible, implying a lower asymptotic number of infected.
This fact can clearly be seen through figure 8, where we plot the asymptotic value of susceptible as a function of the start time $t_I$.
The downside is the sooner the strategy starts,  the longer it takes to finish. We can see this through the right panel of figure 7, where we plot the duration of the strategy as a function of the start time $t_I$. We clearly see that the sooner the control begins, the longer the strategy will remain. Just as an example, note that when $ t_I = 40 $ days the strategy must last more than $ 10,000 $ days, which represents approximately $ 27 $ years.
Therefore, to set up a good $S_0^*$-Strategy one should choose to start the control on a time $t_I$ that  implies a significant reduction of the asymptotic total number of infected and also a significant flattening of the peak, but avoiding a very long duration.

\begin{figure}
[!htb]
\begin{center}
\includegraphics[width= 14.0cm]{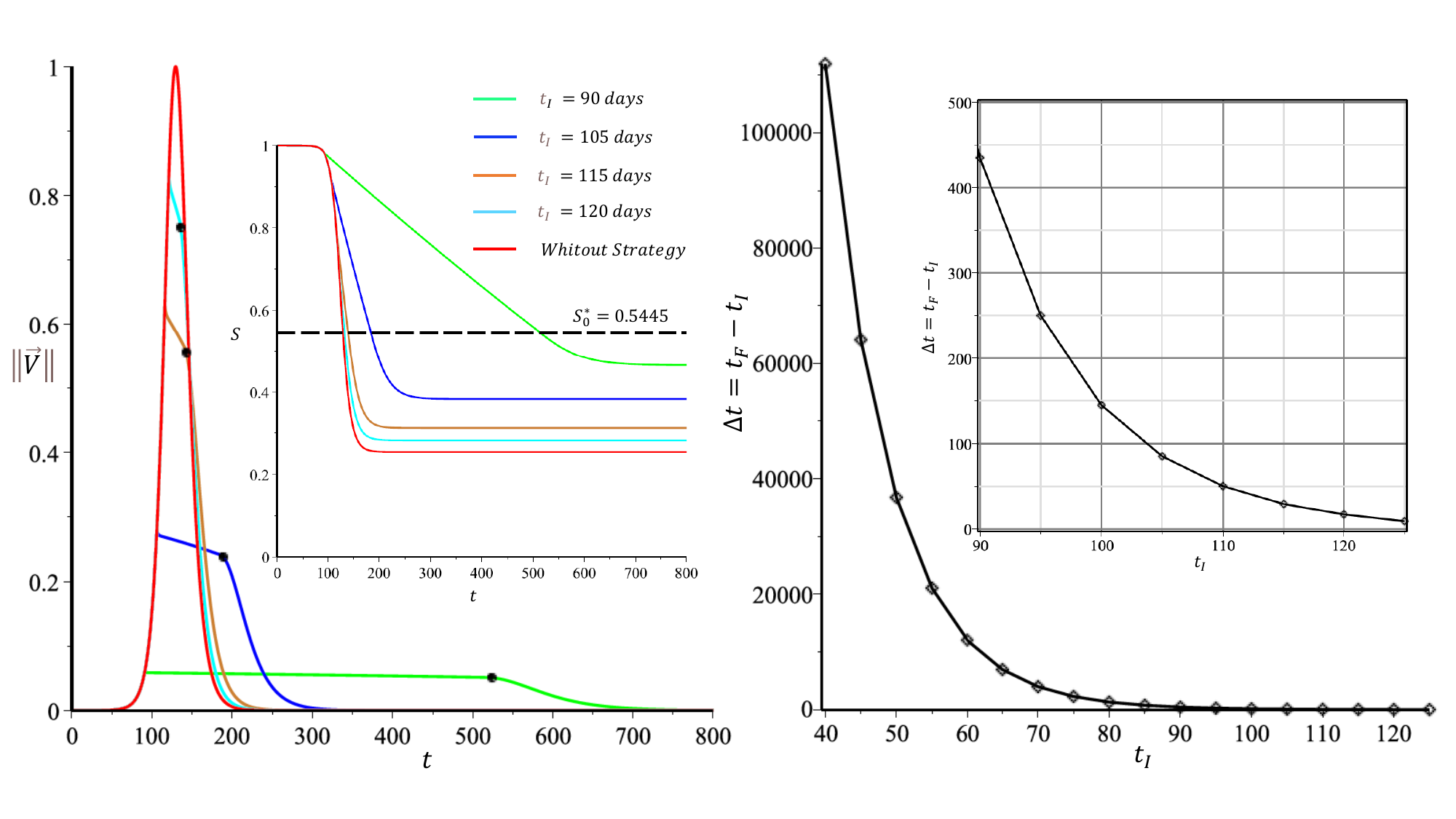}
\end{center}
\vspace*{-10mm}
\caption{$S_0^*$-Strategies with different start time $t_I$}
\end{figure}

\begin{figure}
[!htb]
\begin{center}
\includegraphics[width= 14.0cm]{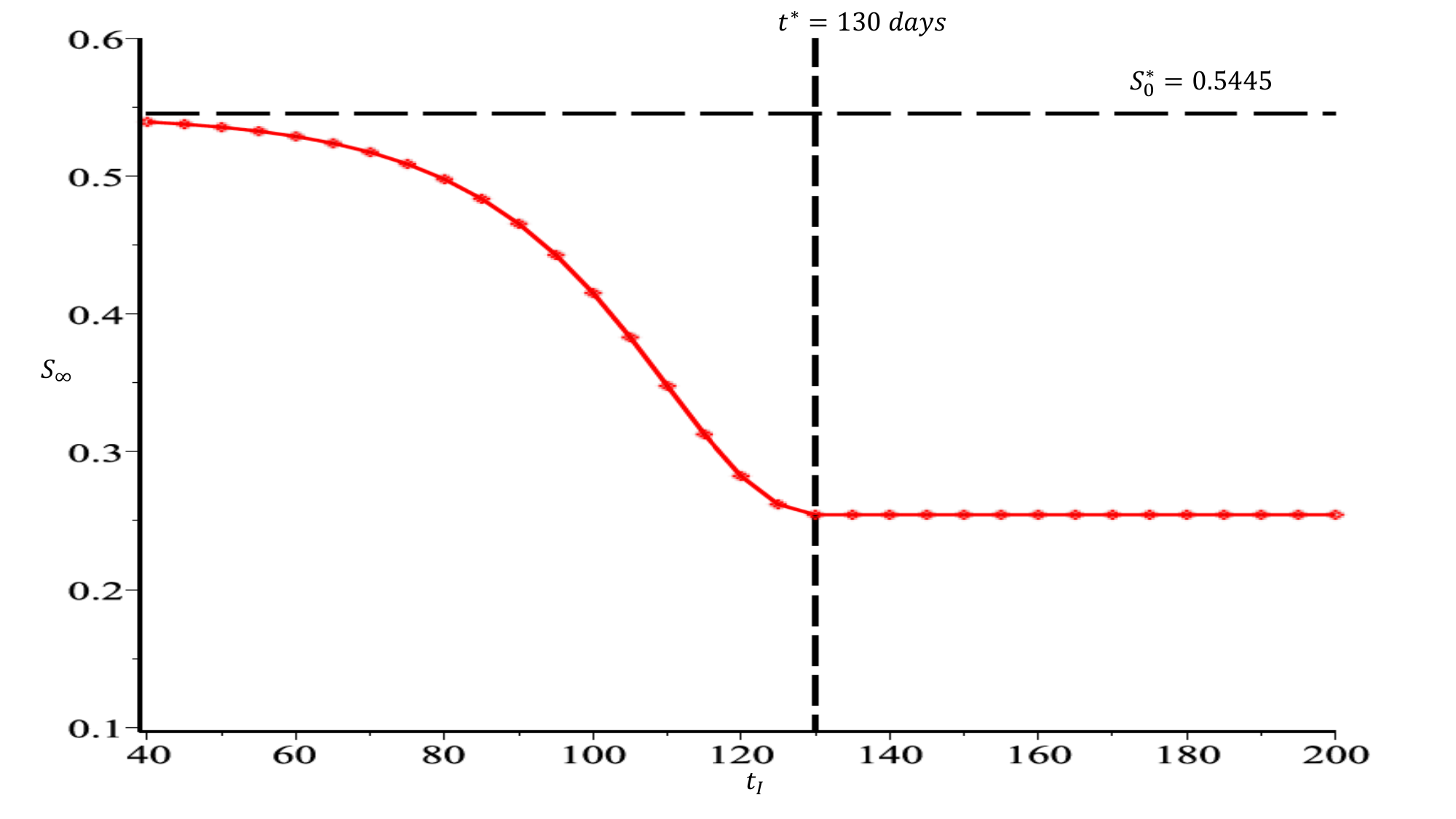}
\end{center}
\vspace*{-10mm}
\caption{Asymptotic Susceptible $S_\infty$ for $S_0^*$-Strategies with different start time $t_I$}
\end{figure}

Finally, we show in figure 9 the relation between the duration $ t_F-t_I $, the asymptotic value $ S_\infty $ and the period $ \Delta $ of $ S_0^* $ -Strategies that start at the same time $ t_I $. In the larger panel on the right, we see that the asymptotic value of susceptible $S_\infty$ is growing and
seems to converge to $ {S_0^*} $ as $ \Delta $ increases.
In the left panel we see that the duration of the strategy increases non-linearly with the increase of the period $ \Delta $.
In the smaller panel on the left we have a zoom of the larger panel and we can see that for $ \Delta = 80 $ days we have a duration of approximately $ t_F-t_I = 1 $ year. In the smaller panel on the right we show the corresponding $ S^*_0 $-strategy. We see that, although the strategy lasts almost one year, after only 100 days in which the value of $ q (t) $ is slightly greater than $ 0.6 $ it jumps to a value greater than $ 0.8 $, that is, after a little more than three months, the control strategy needs to reduce the value of $ R_0 $ by less than $ 20 \% $ of its natural value. Moreover, the value $ 1-S_\infty \approx 0.49 $ corresponds to a value close to $ 1-{S_0^*} \approx 0.46 $, which represents the maximum possible reduction in the asymptotic number of the total infected.

\begin{figure}
[!htb]
\begin{center}
\includegraphics[width= 14.0cm]{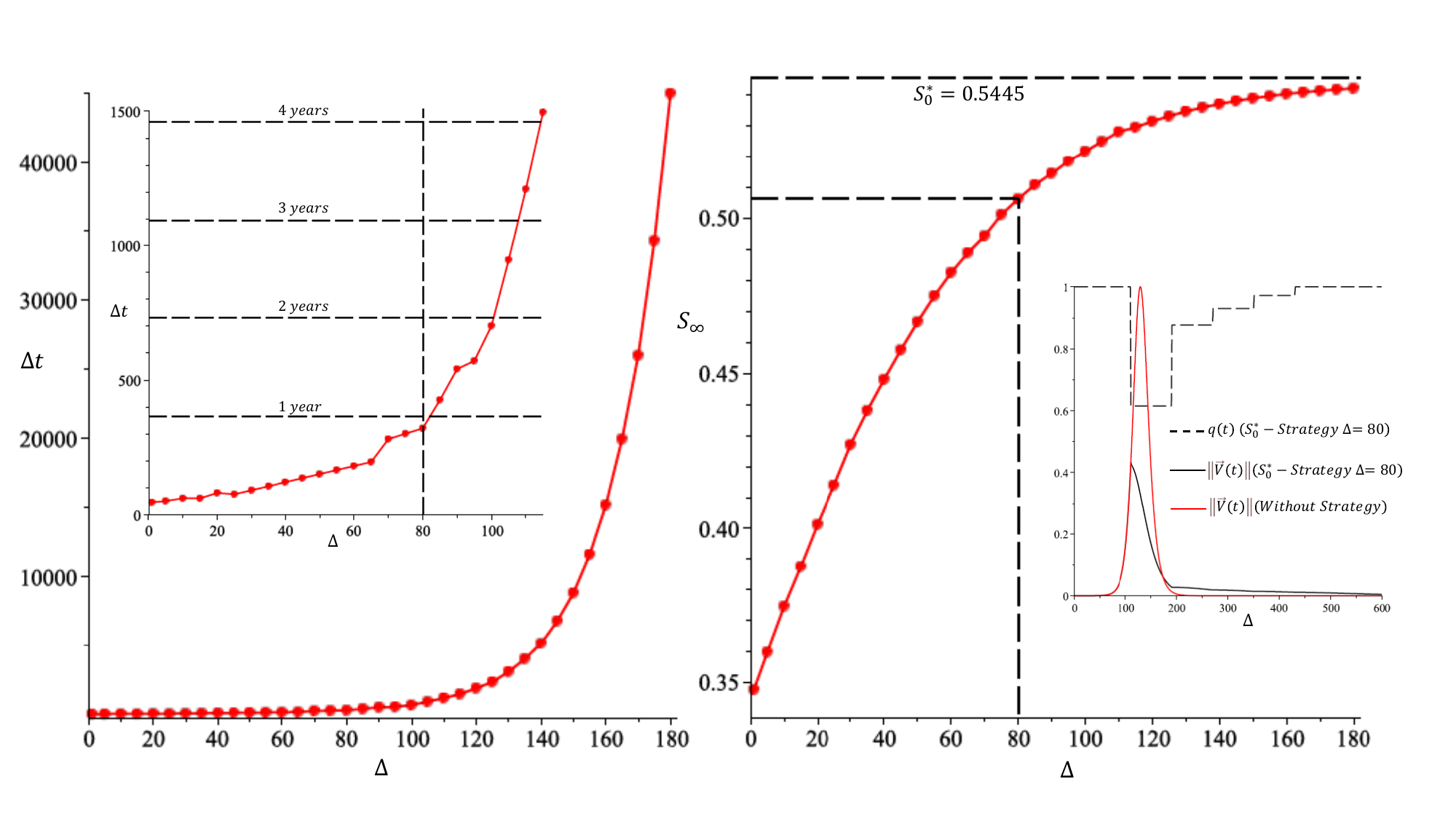}
\end{center}
\vspace*{-10mm}
\caption{$S_0^*$-Strategies with different period $\Delta$}
\end{figure}

The great advantage of $S^*_0$-strategies is to allow a return to the natural situation after a finite time, always maintaining the condition of non-outbreak. However, an important problem relates with the general problem of any regulated strategies: to built a good monitoring system to have reliable information about the current state of the model variables. 
For $S^*_0$-Strategies, it is necessary a good estimation on the value of susceptible throughout the evolution and application of the control strategy.

\section{The SEIAQR Model for COVID-19: A Control Strategy based on Isolation of Infected Individuals}

Uniform control strategies are based on the idea to change in the same proportion the values of  reproduction numbers $ R_i $ ($i=1,\ldots,N$), which are associated to the different infectious stages $ I_i $  in the SEIR model given in (\ref{seir}). An {\it Isolation Strategy (or Quarantine)} can be seen as an extreme case of control, where for a certain group of infected individuals is imposed a zero reproduction number, that is, it is prevented to these individuals to spread the disease.
In this sense, an isolation strategy may be modeled considering three infectious stages: $ I_1 = I $ (Symptomatic Infected), $ I_2 = A $ (Asymptomatic Infected) and $ I_3 = Q $ (Isolated Infected - Quarantine), where the respective reproduction numbers are written as $ R_1 = R_0 $, 
$ R_2 = \xi R_0 $ and $ R_3 = 0 $. As in the case for a SEIAR model one considers equal infection times $ \gamma_1 = \gamma_2 = \gamma_3 = \gamma $,
thus, directly from the model for equal infection times in (\ref{seiarap}), one can write the following model called SEIAQR:
\begin{equation}
\label{seiaqr}
\begin{array}{c}
{\dot S}=-\beta (I+\xi A) S,\;
{\dot E}=\beta (I+\xi A) S-\sigma E,\\ \\
{\dot I}=(1-\zeta_I)(1-\chi)\sigma E-\gamma I,\;
{\dot A}=(1-\zeta_A)\chi\sigma E-\gamma A,\\ \\
{\dot Q}=\left(\zeta_I(1-\chi)+\zeta_A\chi\right)\sigma E-\gamma Q,\;
{\dot R}=\gamma(I+A+Q),
\end{array}
\end{equation}
where $\beta=\gamma R_0$. The new variable and additional probabilities with respect to SEIAR model in (\ref{seiar}) are shown in table 3.
\vspace*{2mm}
\begin{center}
\begin{tabular}{c}
Table 3 \\
\hline\hline 
Additional SEIAQR Model Variable \\
\hline\hline \\
\begin{tabular}{cc}
$Q$ & Population of Isolated Infected Individuals (Quarantine)
\end{tabular} \\ \\
\hline\hline
Additional SEIAQR Model Parameters\\ \hline\hline \\
\begin{tabular}{cc}
$\zeta_I$ &Probability of Isolation of Symptomatic Infected\\ \\
$\zeta_A$ &Probability of Isolation of Asymptomatic Infected\\
\end{tabular} \\ \\
\hline\hline
\end{tabular}
\end{center}

Similarly to a SEIAR model, the surface $ E = I = A = 0 $ is an invariant surface and the column vector $ {\vec V} = [E, I, A] ^ {t} $  satisfies the following equation:
\begin{equation}
\frac{d{\vec V}}{dt}={\bf L}{\vec V},\;\;{\bf L}=\left(
\begin{array}{ccc}
-\sigma & \beta S & \xi\beta S\\ & &\\
(1-\zeta_I)(1-\chi) \sigma & -\gamma & 0 \\ & &\\
(1-\zeta_A)\chi\sigma & 0 & -\gamma 
\end{array}\right).
\end{equation}
The critical value of susceptible ${S^*}$ in (\ref{SLap}) written as
\begin{equation}
{S^*}=\frac{1}{R_0}\frac{1}{\chi\xi(1-\zeta_A)+(1-\chi)(1-\zeta_I)}.
\end{equation}
The condition for a non outbreak of the epidemic at time $t=0$ when $S(0)=1$  is given by
\begin{equation}
\label{SLseiaqr}
{S^*}\geq 1\;\;\Rightarrow\;\;\chi\xi(1-\zeta_A)+(1-\chi)(1-\zeta_I)\leq\frac{1}{R_0}.
\end{equation} 
Assuming that we can isolate all symptomatic infected, that is, considering $ \zeta_I = 1 $, we have
\begin{equation}
\chi\xi(1-\zeta_A)\leq\frac{1}{R_0}\;\;\Rightarrow\;\;\zeta_A\geq 1-\frac{1}{\chi\xi R_0}.
\end{equation}
Finally, using the parameters of COVID-19 (shown in table 2) we obtain
\begin{equation}
{S^*}\geq 1\;\;{\rm e}\;\;\zeta_I=1\;{\rm (COVID-19)}\;\;\Rightarrow\;\;\zeta_A\geq 0.2969134503.
\end{equation} 
This result shows that we need to isolate at least $ 30 \% $ from asymptomatic infected individuals. This is a huge task, taking into account that $ 86 \% $ of the infected are asymptomatic. This implies the necessity to test an important proportion of the total population, or to have a good system to trace the contacts made by an identified infected.

We must note that a non-outbreak of the epidemic based on isolation probabilities given in (\ref{SLseiaqr}) can be difficult to achieve in many practical situations, which should not necessarily imply the abandonment of  isolation strategies (quarantine).
Indeed, supposing the impossibility of achieving a non-outbreak quarantine, we must combine outbreak quarantines with uniform strategies settled to control the basic reproduction rate $ R_0 $. The advantage is to have a value for $ {S^*} $ in (\ref{SLseiaqr}) (defined in the SEIAQR model) greater than the value of $ {S^*} $ in (\ref{SLSH}) (defined in the SEIAR model), which leads to a control strategy (for the rest of the non-isolated population) with softer reduction in the epidemic reproduction rate $R_0$. 

As an example, let us consider that all symptomatic infected are isolated, then, according to equation (\ref{SLseiaqr}), we have ${S^*}=0.7031$ (the value without  quarantine is ${S^*}=0.5445$). Therefore, as a consequence of proposition 3, the maximum value for asymptotic susceptible $S_\infty$ for finite strategies is greater than in the absence o quarantine and, 
for any finite control strategy, there will be a lower total number of  infected $1-S_{\infty}$.

\section{Conclusions}

This article sought to define concepts and establish some properties that help to  analyze control strategies applied with the aim of mitigating the effects of an epidemic.

The uniform strategies are based on the idea to reduce in equal proportions the basic reproduction number associated to the entire population.The main characteristic of this kind of control strategy is to prescribe some general rules without discrimination with regard to any special group of infected. Moreover, these rules  must be designed in order to reduce the basic reproduction rate of the epidemic. As example we can cite: compulsory use of masks; general sanitary measures such as washing hands, sanitizing contact surfaces and so on; to establish rules to keep the distance between individuals and to diminish as possible the social contact among them. The extreme possibility in the latter case is to establish a complete blockade (lock-down) to prevent the movement of the entire population considered.

It is quite clear, from the Fundamental Principle for Finite Strategies (proposition 4), that there is a maximum limit to mitigate the effects of an epidemic using only uniform strategies on the entire population,
because, whenever the control strategy stops and the natural epidemic evolution returns, it is obliged to ``deliver susceptible to the disease agent (virus)''. This fact does not imply that such kind of strategy must be abandoned, because we have showed that, in spite the limitation implied by proposition 4, strong reduction in the asymptotic number of total infected can be achieved with not so long time duration for the control strategy.   
Also, a positive fact for finite uniform strategies is the possibility to control the peak intensity of epidemic: the constant finite strategies break the singular peak of epidemic in two smaller peaks and the $S^*_0$-Strategies put the epidemic in a non-outbreak evolution. These properties  
are very interesting because they lead to a lower number of  simultaneous new infected, which is  good to not stress the healthcare system that need to treat all critical ill individuals.

Another kind of control strategy is based on the isolation of infected individuals (quarantine), which means to set to zero their corresponding reproduction numbers. 
Indeed, quarantine is a kind of very ancient control strategy that relies on the possibility to identify and to track infected individuals and the possibility to apply a good quarantine strategy is based on the reliance and efficiency of a good system of epidemiological vigilance. The major challenge for an epidemic like COVID-19 is the evasion of asymptomatic (or symptomatic with mild symptoms) from the tracking system of a given system of epidemiological vigilance.
From a mathematical (or physical) point of view, the main advantages of quarantine is to incorporate a new stage (variable) in the model (\ref{seir}) without to change the reproduction rates of the other stages. This leads to a epidemic evolution that does not alter the critical time $t^*$ when $S(t^*)={S^*}$, although the value of ${S^*}$ itself may be increased. This allows to softening the effect due to propositions 2 and 4 without to set up a control strategy that leads to a higher time duration than the time duration of the epidemic cycle.
In this sense, the only issue with respect to quarantine is the enforcement to keep the functioning of the epidemiological vigilance while the epidemics evolves.
Provided the impossibility to built a quarantine strategy that satisfies a non-outbreak condition, the better is to combine it with a finite  uniform strategy. In this case, the existence of quarantine strategy allows to built a uniform control strategy with less severe restrictions for the non isolated individuals.

One last observation that we consider important is on the idea of herd immunity, more specifically with the idea of a value for the proportion of susceptible  from which the epidemic stops its spread. As shown in \cite{ref36,ref37} this value depends on the basic reproduction rate $ R_0 $ and has been calculated for different countries.
The point here is that this value should correspond to the asymptotic value $ S_\infty $ associated with the model defined in (\ref{seir}), however this value is not unique, since it corresponds just to one of the many possible fixed points of the model. This is quite clear when we see how the application of different finite control strategies implies a different value for 
$ S_\infty $, which is directly related to the fact that these asymptotic fixed points depends on the initial condition. Indeed, the only general conclusion we can draw is about the maximum possible value for $ S_\infty $.

\section{Acknowledgments}

This research did not receive any specific grant from funding agencies in the public, commercial, or not-for-profit sectors.

\begin{appendix}

\section{The Critical Value of Epidemic} 

The invariant surface $  E =J= 0 $ of the differential system (\ref{seired})  is constituted by an infinity number of fixed points, each of them characterized by a different value of $S$. 
If we define the column vector $ {\vec v} = [E, J] ^ {t} $, then we can write
\begin{equation}
\frac{d{\vec v}}{dt}={\bf l}{\vec v},\;\; {\bf l}(S)=\left(
\begin{array}{cc}
-\sigma & \beta S \\
\bar\sigma & -\gamma
\end{array}\right).
\end{equation}

Let us consider the vector field $ {\bf l}{\vec v} = [\dot {E}, \dot {J}]^t $, defined in the bi-dimensional $ (E, J)\in I\hspace*{-1.5mm}R^{2}$ space, for a fixed value of $ S $. The sign analysis of $ \dot {E} $ and $ \dot {I} $ leads to:
\begin{equation}
\label{sinalEJ}
{\rm sgn}(\dot{E})=\left\{\begin{tabular}{cc}
$-1$ &  $\displaystyle \frac{J}{E}<\frac{\sigma}{\beta S}$\\ & \\
$0$ &  $\displaystyle \frac{J}{E}=\frac{\sigma}{\beta S}$ \\ & \\
$+1$ &  $\displaystyle \frac{J}{E}>\frac{\sigma}{\beta S}$ \end{tabular}\right.\; ,\;
{\rm sgn}(\dot{J})=\left\{\begin{tabular}{cc}
$+1$ &  $\displaystyle \frac{J}{E}<\frac{\bar\sigma}{\gamma}$\\ & \\
$0$ &  $\displaystyle \frac{J}{E}=\frac{\bar\sigma}{\gamma}$ \\ & \\
$-1$ &  $\displaystyle \frac{J}{E}>\frac{\bar\sigma}{\gamma}$ \end{tabular}\right. .
\end{equation}
The definition of the following two angles:
\begin{equation}
\tan(\theta_E)=\frac{\sigma}{\beta S},\;\;\tan(\theta_J)=\frac{\bar\sigma}{\gamma},
\end{equation}
allows to represent, through the sign analysis made in (\ref{sinalEJ}), the direction of the vector field $ {\bf l}{\vec v}$ in the $(E,J)$ plane.
This representation is shown in the figure A.10, respectively for the cases $ \theta_E <\theta_J $ and $ \theta_E> \theta_J $. 

\begin{figure}
[!htb]
\begin{center}
\includegraphics[width= 13.0cm]{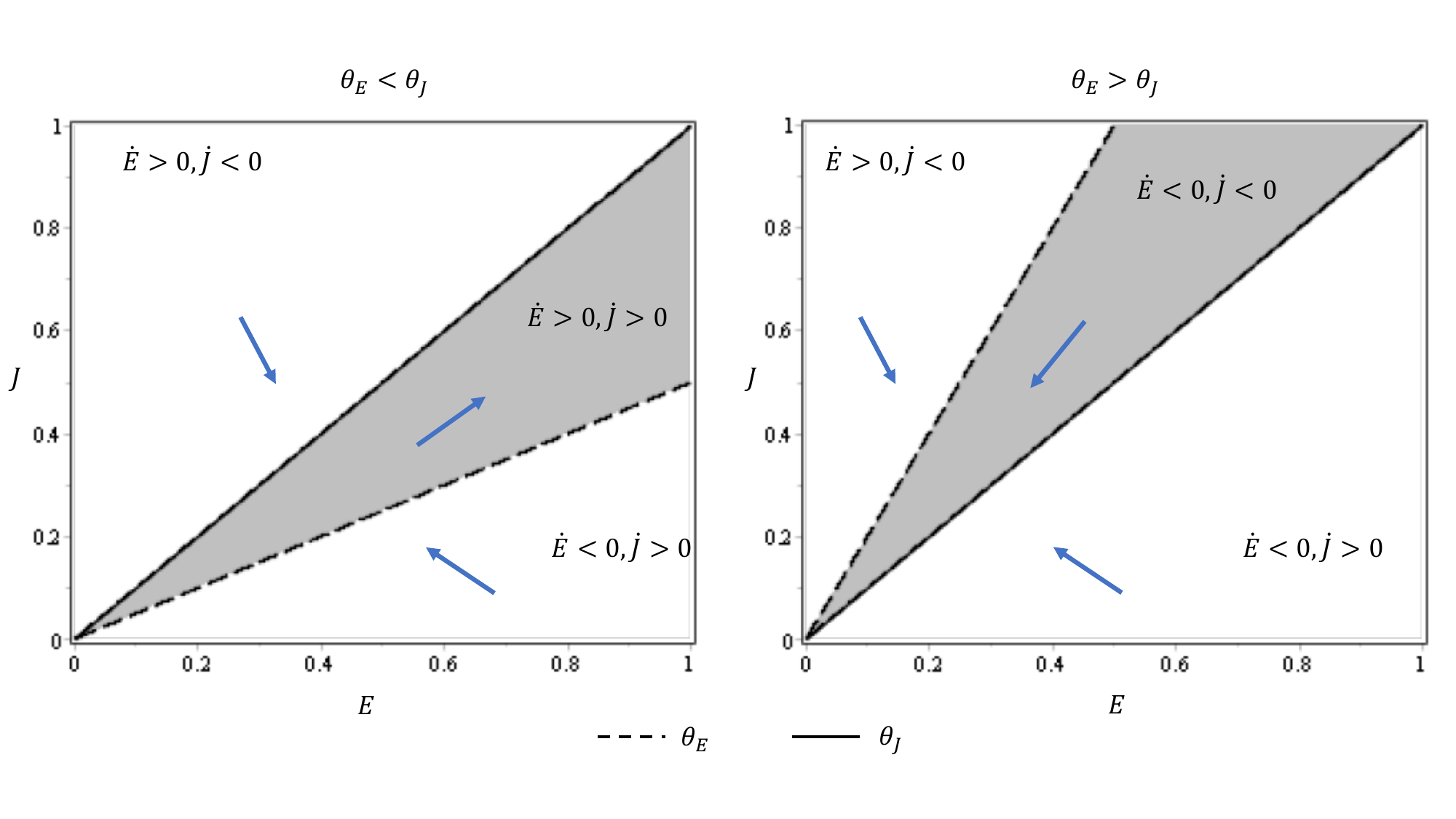}
\end{center}
\vspace*{-10mm}
\caption{Topological transition in the field ${\bf l}{\vec v}$.}
\end{figure}

We clearly see that  $ \theta_E <\theta_J $ leads to $ S> \gamma\sigma / \beta\bar\sigma $ and $ \theta_E> \theta_J $ implies in $ S <\gamma\sigma / \beta\bar\sigma $.
Also, we see that these two inequalities represent two different topological configurations in the vector field ${\bf l}{\vec v}$ directions and 
the transition between them happens when $ S= S_{c} $:
\begin{equation}
\label{transicao}
S_{c}=\frac{\gamma\sigma}{\beta\bar\sigma}=\frac{1}{\sum_{i=1}^Nx_ie_i}\frac{1}{R_0}={S^*},
\end{equation}
where $\tan (\theta_E) =\tan (\theta_J)$.
It is worth to observe that for $ S> {S^*} $ the vector field ${\bf l}{\vec v}$ implies trajectories that tend to depart from the origin of  $ (E, J)$ plane and for $ S <{S^*} $ they tend to approach this origin. Therefore, the time $t^*>0$ such that $ S(t^*)={S^*} $ may represent a turning point in the epidemic evolution.

An initial condition for a given time instant $ t_0> 0 $ such that
\begin{equation}
\label{condinicial}
S(t_0)>{S^*},\;\tan(\theta_E)\leq \frac{J(t_0)}{E(t_0)}\leq\tan(\theta_J),
\end{equation}
corresponds to a point $ (E (t_0), J (t_0)) $ in the gray region of the left panel of figure A10. The straight lines in this figure: the dashed line and the solid line, respectively represent the region where $ \dot E = 0 $ and $ \dot J = 0 $. While $ S (t)> {S^*} $ for $ t> t_0 $, the vector field $ [\dot E, \dot J]^t $ will be pointing into the gray region. This property stems from the fact that on the dashed line we have $ \dot E = 0 $ and $ \dot J> 0 $ and on the solid line $ \dot J = 0 $ and $ \dot E> 0 $. 

In proposition 2 we have shown that $ \lim_ {t \rightarrow \infty} S (t) < {S^*} $, so there is a time $ t^*> t_0 $ such that $ S (t^*) = {S^*} $. Therefore, for the time interval $ t_0 \leq t <t^* $ a solution that starts in an initial condition that satisfies (\ref{condinicial}) cannot leave the gray region (left panel of figure A10), because, as shown in the previous paragraph, the vector field $ [\dot E, \dot J] ^ t $  at the borders of this gray region is always pointing to inside it. Thus, we can affirm that in the interval $ t_0 <t <t^* $ we have $ {\dot E} (t)> 0 $ and $ {\dot J} (t)> 0 $;  and  for $ t = t^* $ we have $ {\dot E} (t^*) = {\dot J} (t^*) = 0 $.

Considering the functions $ \tan \left(\theta (t) \right) = E (t) / J (t) $, $ \tan \left (\theta_E (t) \right) = \sigma / \beta S (t) $ and the angle $ \tan (\theta_L) = \sigma / \beta {S^*} = {\bar \sigma} / \gamma $, we can show that for a time $ t_1 = t^* + \delta t $ it holds the following expansions in power series of $ \delta t $:
$$\tan\left(\theta(t_1)\right)=\tan\left(\theta_L\right)+A{\delta t}^2+O\left(\delta t^3\right),$$
$$\tan\left(\theta_E(t_1)\right)=\tan\left(\theta_L\right)+B\delta t+O\left(\delta t^2\right),$$
where
$$A=-\beta\frac{J^2(t^*){\dot S}(t^*)}{E^2(t^*)},\;B=-\frac{\sigma}{\beta}\frac{{\dot S}(t^*)}{{S^*}^2}.$$
Due to the fact that $ {\dot S} (t^*) = - \beta J (t^*) {S^*} <0 $ we have $ A> 0 $, $ B> 0 $ and, consequently, the angles $ \theta (t) $ and $ \theta_E (t) $ are increasing functions for a sufficiently small range of $ \delta t $. Furthermore, there is a value $ \epsilon> 0 $ such that for 
$ 0 <\delta t <\epsilon $ we have
$$\theta_J<\theta(t_1)<\theta_E(t_1).$$
This means that the trajectory of the solution considered so far will be in the gray region of the right panel of figure A10.
Moreover, $ S (t) <{S^*} $ for all $ t> t^* $ and, analogous as shown for the gray region of the left panel of figure A10, we can show that the trajectory of the considered solution cannot leave the gray region in the right panel of the figure. This fact occurs because now in the region where $ \dot E = 0 $ (dashed line) we have $ \dot J <0 $ and in the region where $ \dot J = 0 $ we have $ \dot E <0 $, showing that at the borders of the gray region the vector field $ [\dot E, \dot J] ^ t $ points to inside it.
Finally, we can conclude that $ {\dot E} (t) <0 $ and $ {\dot J} (t) <0 $ for all time $ t> t^* $.

We can summarize the results obtained in the following proposition:

\begin{proposicao}
For any solution with an initial condition that satisfies the relations in (\ref{condinicial}), there is a time $ t^*> t_0 $ such that
\begin{itemize}
\item ${\dot E}(t)>0$ and ${\dot J}>0$ for $t_0<t<t^*$.
\item ${\dot E}(t)=0$ and ${\dot J}=0$ for $t=t^*$.
\item ${\dot E}(t)<0$ and ${\dot J}<0$ for $t^*<t$.
\end{itemize}
\end{proposicao}

We define an initial condition at an instant time $ t_0 $: $ S (t_0) $, $ E (t_0) $ and $ J (t_0) $, as an {\it Initial Outbreak Condition}, if there is a time $ t_1 \geq t_0 $ such that
$$ S (t_1)> {S^*}, \; \; \frac{\sigma} {\beta S (t_1)} = \tan \left (\theta_E (t_1) \right) \leq \frac {J (t_1)} {E (t_1)} \leq \tan \left (\theta_J \right) = \frac{{\bar \sigma}} {\gamma}. $$
We can see that an initial outbreak condition necessarily implies that $ S (t_0)> {S^*} $. This situation is represented in the left panel of figure A10 and we can conclude that all initial conditions in the gray region are initial outbreak conditions with $ t_1 = t_0 $. On the other hand, for initial conditions of outbreak outside the gray region, there is a time $ t_1 $ such that the trajectory of the corresponding solution will enter the gray region before $ S (t) $ becomes less than $ {S^*} $. Therefore,
it is easy to conclude that for all solution with initial outbreak conditions, there is a time $ t_2\geq t_1 $ such that $ \dot E (t_2)> 0 $ and $ \dot J (t_2)> 0 $.

\end{appendix}


\begin{thebibliography}{99}
\bibitem{ref38} In reference \cite{ref16} the value of $R_0$ estimated in China range from $2.24$ to $3.58$. In reference \cite{ref18} the value of $R_0$ estimated for the state of Bahia-Brazil range from $2.25$ to $3.56$.  In reference \cite{ref39} the value of $R_0$ estimated for the city of Brasília-Brazil range from $2.5$ to $3.5$. We pick up  $R_0=3.0$ to represent the value for the basic reproduction rate for Brasília.
\bibitem{ref42} By asymptomatic we understand all  individuals with no symptoms or mild symptoms that not allow to distinguish COVID-19 from other  respiratory diseases and can not easily be tracked unless specific laboratory tests are performed (see references \cite{ref15}, \cite{ref28}, \cite{ref34} and \cite{ref35}).
\bibitem{ref43} This system of ODE's  is obtained from the model presented in reference \cite{ref33} in the non endemic condition, that is, without to consider the birth an mortality rates. Also we do not take into account the amelioration rates that appears in this model.
\end{thebibliography}
\end{document}